\documentclass[copyedit]{informs1}  % for review, blinded AND (double)-spaced for copyediting

\usepackage{endnotes}
\usepackage{enumerate}
\usepackage{rotating}
\usepackage{subfigure}
\usepackage{graphicx}
\usepackage{multirow}
\usepackage{hhline}
\newcommand{\mb}[1]{\mbox{\boldmath $#1$}}
\newcommand{\mbs}[1]{{\mbox{\boldmath \scriptsize{$#1$}}}}
\newcommand{\mbss}[1]{{\mbox{\boldmath \tiny{$#1$}}}}

\usepackage{setspace,graphics,booktabs,color}
\usepackage{multirow}

\usepackage{natbib}
 \bibpunct[, ]{(}{)}{,}{a}{}{,}%
 \def\bibfont{\small}%
\makeatletter

\newcommand{\Rmnum}[1]{\expandafter\@slowromancap\romannumeral #1@}
\makeatother

\newcommand{\bx}{\mbox{\boldmath $x$}}
\newcommand{\bu}{\mbox{\boldmath $u$}}
\newcommand{\by}{\mbox{\boldmath $y$}}
\newcommand{\bv}{\mbox{\boldmath $v$}}
\newcommand{\bb}{\mbox{\boldmath $b$}}

\newcommand{\bSigma}{\mbox{\boldmath $\Sigma$}}

\newcommand{\bz}{\mbox{\boldmath $z$}}
\newcommand{\bd}{\mbox{\boldmath $d$}}
\newcommand{\bmu}{\mbox{\boldmath $\mu$}}
\newcommand{\bepsilon}{\mbox{\boldmath $\epsilon$}}
\newcommand{\bq}{\mbox{\boldmath $q$}}

\TheoremsNumberedThrough     % Preferred (Theorem 1, Lemma 1, Theorem 2)
\EquationsNumberedThrough    % Default: (1), (2), ...
\begin{document}
%%%%%%%%%%%%%%%%

\RUNAUTHOR{}
%
%% Title or shortened title suitable for running heads. Sample:
%% \RUNTITLE{Bundling Information Goods of Decreasing Value}
%% Enter the (shortened) title:
\RUNTITLE{}

% \TITLE{Bundling Information Goods of Decreasing Value}
\TITLE{\Large Analysis of Discrete Choice Models: A Welfare-Based Framework}

\ARTICLEAUTHORS{\AUTHOR {Guiyun Feng} \AFF{Department of Industrial
and Systems Engineering, University of Minnesota, Minneapolis, MN
55455, {fengx421@umn.edu}}

\AUTHOR {Xiaobo Li} \AFF{Department of Industrial and Systems
Engineering, University of Minnesota, Minneapolis, MN 55455,
{lixx3195@umn.edu}}

\AUTHOR {Zizhuo Wang} \AFF{Department of Industrial and Systems
Engineering, University of Minnesota, Minneapolis, MN 55455,
{zwang@umn.edu}}}
%\ABSTRACT{Observing that the consumer choice vector in the random utility model can be derived as the gradient of the corresponding expected utility function, we provide a general framework to obtain consumer choice by axiomatizing the expected utility function. We assume that the expected utility functions are monotone, convex and translation invariant, all of which are satisfied in random utility model. We find that this axiomatic framework is in fact equivalent to the other two existing framework: representative agent model and semi-parametric model. We note that the equivalence of the latter two models is not known in the literature. While in $n\geq 3$ cases, the representative agent model is known (for example, in \cite{hofbauer}) to be strictly more general to random utility model, we prove that in the case of $2$ products, they are equivalent. Lastly, we study substitutability and complementarity properties of the choice model. We show that our generalizations of the random utility model can model complementarity between products, while in random utility models, all products are substitutable to each other.    }
% \KEYWORDS{}

\ABSTRACT{Based on the observation that many existing discrete
choice models admit a  welfare function of utilities whose gradient
gives the choice probability vector, we propose a new representation
of discrete choice model which we call the \textit{welfare-based
choice model}. The welfare-based choice model is meaningful on its
own by providing a new way of constructing choice models. More
importantly, it provides great analysis convenience for establishing
connections among existing choice models. We prove by
using convex analysis theory, that the welfare-based choice model is
equivalent to the representative agent choice model and the
semi-parametric choice model, establishing the equivalence of the
latter two. We show that these three models are all strictly more
general than the random utility model, while when there are only two
alternatives, those four models are equivalent. In particular, we
show that the distinction between the welfare-based choice model and
the random utility model lies in the requirement of the higher-order
derivatives of the welfare function. We then define a new concept in
choice models: substitutability/complementarity between
alternatives. We show that the random utility model only allows
substitutability between different alternatives; while the
welfare-based choice model allows flexible
substitutability/complementarity patterns. We argue that such
flexibility could be desirable in capturing certain practical choice
patterns and expanding the scope of discrete choice models. Examples
are given of new choice models proposed under our framework. }
\HISTORY{This version: \today} \maketitle

\section{Introduction}
\label{sec:introduction}

In this paper, we study the discrete choice models. Discrete choice
models are used to model choices made by people among a finite set
of alternatives. For example, they are used to examine which product
to purchase for a consumer, which mode of transportation to take for
a passenger, among many other choice scenarios people face everyday.
In the past few decades, discrete choice models have attracted great
interest in the economics, marketing, operations research and
management science communities. Specifically, such models have been
viewed as the behavioral foundation in many operational
decision-making problems, such as transportation planning,
assortment optimization, multiproduct pricing, etc.

In the past few decades, researchers have proposed a variety of
discrete choice models (see \citealt{andersonbook} and
\citealt{Ben-Akiva_Lerman}). Among them, the most popular one is the
random utility model, in which a utility is assigned to each
alternative. In the random utility model, the utility is composed of
a deterministic part and a random part. Each individual then chooses
the alternative with the highest utility, given the realization of
the random part. Different choice models arise when different
distributions for the random part are used. Some examples of random
utility model can be found in \cite{McFadden1974,McFadden_nested}
and \cite{Daganzo}. Another popular choice model is the
representative agent model, in which a representative agent makes
the choice on behalf of the population. In the representative agent
model, there is again a utility associated with each alternative,
and the representative agent  maximizes a weighted utility of the
choice (which is a vector of proportions for each alternative) plus
a regularization term, which typically encourages diversification of
the choice (\citealt{anderson_representative}). More recently, a
class of semi-parametric models has been proposed (see
\citealt{natarajan2009persistency}). This model is similar to the
random utility model. However, instead of specifying a single
distribution for the random utility, a set of distributions is
considered. Then they choose one extreme distribution in that set
to determine the choice probabilities. There are other choice models
based on the dynamics of choice decisions or other non-parametric
ideas. We will provide a more detailed review of these models in
Section \ref{sec:review}.

Although these models have all provided excellent explanations, both
theoretically and empirically, for how people make choices in
practice, some gaps in the literature still exist regarding the
relations between those popular choice models. In particular, the
following questions are not answered in the prior literature:
\begin{enumerate}
\item What is the relation between the representative agent
model and the semi-parametric model? It has been shown that for many
special cases, the semi-parametric model can be represented as a
representative agent model. However, it is unknown whether this is
generally true.
\item It is known that both the representative agent model and the semi-parametric model are more
general than the random utility model. What exactly is the
distinction between these models?
\item What choice pattern is restricted in the  random utility model? Can we easily
construct choice models that relax those restrictions?
\end{enumerate}

In this paper, we present precise answers to the above questions. To
answer those questions, we propose a new class of choice models,
which we call the {\it welfare-based choice model}. The
welfare-based choice model is based on the observation that many
existing choice models take the form of mapping a utility vector to
a probability vector and admit a welfare function of the utilities
whose gradient gives the choice probability vector. Therefore, by
directly proposing desirable conditions on the welfare functions, we
define the class of welfare-based choice models. We show that the
welfare-based choice model is not only meaningful on its own, but
also provides great analysis convenience for establishing
connections between existing choice models.

First, we show that by using the welfare-based choice model as an
intermediate model, the classes of choice models defined by: 1) the
welfare-based choice model, 2) the representative agent model and
3) the semi-parametric model, are the same. More precisely, under
mild regularity assumptions, given any of the following
three: a  choice welfare function (which defines a welfare-based
choice model), a regularization function (which defines a
representative agent model) or a distribution set (which defines a
semi-parametric model), one can construct the other two to define
exactly the same choice model. This means that the class of
representative agent models and the class of semi-parametric models
are equivalent to each other, which is somewhat surprising because
they seem to have very different origins. In addition, our proof of
the equivalence of these three models is constructive, therefore,
it gives methods to convert one model to another in an explicit way,
potentially alleviating the need to construct correspondences in a
case by case manner as is done in the current research.

Second, we study the relation between the above three models and the
random utility model. We show that when there are only two
alternatives, the random utility model is equivalent to the above
three models. We also demonstrate that this is not true in general
if there are more than two alternatives, in which case the above
three models strictly subsume the random utility model. In
particular, we point out the exact distinction between these three
models and the random utility model, which lies in the higher-order
derivatives of the  choice function. Our result gives  precise
relations among those models.

Finally, by examining the difference between the welfare-based
choice model and the random utility model, we identify an important
property that is restricted in the random utility model but is
flexible in the other three models. We call the property {\it
substitutability and complementarity of alternatives}. Specifically,
this property examines whether the choice probability of another alternative will
increase or decrease when the utility of
one alternative increases. We show that random utility models only allow
substitutability between alternatives. Although this is natural in
many practical situations, we argue that in certain applications, it
might be appealing to allow some alternatives to exhibit
complementarity in certain range, especially when the utility
is based on scores and certain alternatives share the same
feature (e.g., brand or certain component) on which the score is based. We
derive conditions under which different choice models exhibit substitutable/complementary properties. As
far as we know, we are the first to study such properties in choice
models, and we believe that this study will open new possibilities
in the design of choice models by enlarging its horizon and
capturing more practical choice patterns. In fact, we show a few
examples of  new choice models that allow complementarity among
choices (in a certain range) and explain the practicality of those
models.

It is worth mentioning that the analysis technique used in this
paper is novel and interesting. In particular, we adopt several convex
analysis tools that are not commonly used in the study of choice
models. Such tools enable us to uncover  deep connections
between seemingly unrelated models and are  key to our findings.
We believe that such analysis methods may be of independent interest
in the future study of choice models.

The remainder of this paper is organized as follows: In Section
\ref{sec:review}, we review discrete choice models that are relevant
to our study. In Section \ref{sec:unified_framework}, we propose the
welfare-based choice model and study its relation with other choice
models. In Section \ref{sec:relation_rum}, we study the relation
between the welfare-based choice model and the random utility model.
In Section \ref{sec:substitutability}, we propose the concept of
substitutability and complementarity between choice alternatives and
derive conditions under which each model exhibits such properties.
We discuss the issue of constructing new choice models from existing
ones in Section \ref{sec:construction}. We conclude the paper in
Section \ref{sec:conclusion}.

\textbf{Notations.} Throughout the paper, the following notations
will be used. We use notation $\mathcal R$ to denote the set of real numbers, and $\bar{\mathcal R}=\mathcal{R}\cup\{-\infty,+\infty\}$ to denote the set of extended real numbers. We use $\mathbf{e}$ to denote a vector of all ones, $\mathbf{e}_{i}$ to denote a vector of zeros except 1 at the $i$th entry, and $\mathbf{0}$ to denote a vector of all zeros (the dimension of these vectors will be clear from the context). Also, we write $\bx \geq \by$ to denote a componentwise relationship and  $\Delta_{n-1}$ to denote the $n-1$-dimensional simplex, i.e., $\Delta_{n-1}=\{\bx|\mathbf{e}^{T}\bx=1, \bx\geq \mathbf{0}\}$. In  our discussions, ordinary lowercase letters $x,y,\dots$ denote
scalars, boldfaced lowercase letters $\bx,\by,\dots$ denote
vectors.
%
% First, although there have been extensive
%studies of certain particular discrete choice models, there lacks
%studies about the
%relationships between different choice models. %In the literature,
%there are often studies about the relation between a few particular
%choice models but not about classes of choice models.
%In particular, many choice models discussed above take a form of
%mapping a vector of utilities to a vector of choice probabilities.
%It is unclear whether a choice model originated from one model can
%be represented through another and how;
%\end{enumerate}
%
%In this paper, we answer the above questions. We first propose a
%unified framework that connects several of the existing choice
%models. The framework is based on the observation that many existing
%choice models take the form of mapping a utility vector to a
%probability vector and admit a choice welfare function whose gradient is
%the choice probability vector.

\section{Review of Existing Discrete Choice Models}
\label{sec:review}

In this section, we review several prevailing classes of discrete
choice models that are related to the discussion in this paper.

\subsection{Random Utility Model}
\label{subsec:rum}

Perhaps the most popular class of discrete choice model is the
random utility model (RUM), proposed first by
\cite{thurstone1927law} and later studied in a vast literature in
economics (see \citealt{andersonbook} for a comprehensive review).
In such a model, a random utility is assigned to each of the
alternatives, and an individual will pick the alternative with the
highest realized utility. Here, the randomness could be due to the
lack of information of the alternatives for a particular individual
or to the idiosyncrasies of preferences among a population. As the
output, the random utility model predicts a vector of choice
probabilities among the alternatives, rather than a single
deterministic choice. Mathematically, suppose there are $n$
alternatives denoted by $\mathcal{N}=\{1,2,...,n\}$, then the random
utility model assumes that the utility of alternative $i$ takes the
following form:
\begin{eqnarray}\label{rum}
u_i = \mu_i+ \epsilon_i, \ \ \forall i \in \mathcal{N},
\end{eqnarray}
where $\bmu = (\mu_1,...,\mu_n)$ is the deterministic part of the
utility and $\bepsilon = (\epsilon_1,...,\epsilon_n)$ is the random
part. In the random utility model, it is assumed that the joint
distribution $\theta$ of $\bepsilon = (\epsilon_1,...,\epsilon_n)$
is known. Then the probability that alternative $i$ will be chosen
is (to ensure the following equation is well-defined, we assume
$\theta$ is absolutely continuous, an assumption we make for all
the random utility models we discuss later):
\begin{eqnarray}\label{rum_choice}
q_i(\bmu)=\mathbb{P}_{\bepsilon\sim \theta}\left(i=\underset{k \in
\mathcal{N}}{\textrm{argmax}}\ (\mu_{k}+\epsilon_{k})\right).
\end{eqnarray}

Random utility models can be further classified by the distribution
function of the random components. The most widely used one is the
multinomial logit (MNL) model, first proposed by
\citet{McFadden1974}. The MNL model is derived by assuming that
$(\epsilon_1,...,\epsilon_n)$ follow independent and identically
distributed Gumbel distributions with scale parameter $\eta$. Given
that assumption, the choice probability in (\ref{rum_choice}) can be
further written as follows:
\begin{eqnarray*}
q_{i}^{{\mathrm{mnl}}}(\bmu)=\frac{\textrm{exp}({\mu_{i}/\eta})}{\underset{k\in
\mathcal{N}}\sum \textrm{exp}({\mu_{k}/\eta})}.
\end{eqnarray*}
It can also be computed that the expected utility an individual can
get under the MNL model is:
\begin{eqnarray*}
w^{\mathrm{mnl}}(\bmu) = \mathbb{E}_{\bepsilon\sim\theta}
\left[\max_{i \in \mathcal{N}} \ \ \mu_i + \epsilon_i  \right] =
\eta \log{\left(\sum_{i\in \mathcal{N}} \exp(\mu_i/\eta)\right)}.
\end{eqnarray*}

The existence of closed-form formulae for the MNL model makes it a
very popular choice model. We refer the readers to
\citet{Ben-Akiva_Lerman} \citet{andersonbook} and
\cite{train2009discrete} for more discussions on the properties of
the MNL model. In addition to the MNL model, there are other choices
of the random part in (\ref{rum}) that lead to alternative  choice
models. Some popular ones among them are the probit model (in which
$\bepsilon$ is chosen to be a joint normal distribution, see, e.g.,
\citealt{Daganzo}), the nested logit model (in which $\bepsilon$ is
chosen to be correlated general extreme value distributions, see,
e.g., \citealt{McFadden_nested}), the mixed logit model (where
$\bepsilon$ is chosen to be Gumbel distributions with a correlated
term, see, e.g., \citealt{McFadden_mixed}, and
\citealt{train2009discrete}) and the exponomial choice model (in
which $\bepsilon$ is chosen to be negative exponential
distributions, see \citealt{Alptekinoglu}).

\subsection{Representative Agent Model}
\label{subsec:ram}

Another popular way to model choice is to use a representative agent
model (RAM). In such a model, a representative agent makes a choice
among $n$ alternatives on behalf of the entire population. In
particular, this agent may choose any fractional amount of each
alternative, or equivalently, his choice is a
vector $\bx = (x_1,...,x_n)$ on $\Delta_{n-1}$. To make his
choice, the agent takes into account the expected utility while
preferring some degree of diversification. More precisely, the
representative agent solves an optimization problem as follows:
\begin{eqnarray}\label{representative}
\mbox{maximize}_{\bx\in \Delta_{n-1}} & \bmu^T\bx - V(\bx).
\end{eqnarray}
Here $\bmu = (\mu_1,...,\mu_n)$ is the deterministic utility of each
alternative, which is similar to that in the random utility model.
$V(\bx): \mathcal{R}^{n}\mapsto \mathcal{R}$ is a regularization term that rewards diversification.
Later, we denote the optimal value of (\ref{representative}) by
$w^r(\bmu)$, which is the  utility a representative agent can obtain if the deterministic
utility vector is $\bmu$. Moreover, if for any $\bmu$, there is a
unique solution to (\ref{representative}), then we define
\begin{eqnarray}\label{representativechoice}
\bq^r(\bmu) = \arg\max\left\{\bmu^T\bx - V(\bx) \big{|} \bx \in
\Delta_{n-1}\right\}
\end{eqnarray}
to be the choice probability vector given by the representative
agent model. (To ensure the maximum is attainable, $V(\bx)$ is
required to be lower semi-continuous. We make this assumption for
all the representative agent models in our discussions.)

A recognized close connection exists between the random
utility model and the representative agent model. In
\cite{anderson_representative}, the authors show that the choice
probabilities from an MNL model with parameter $\eta$ can be equally
derived from a representative agent model with $V(\bx) = \eta
\sum_{i=1}^n x_i \log{x_i}$. Or equivalently, we can write
\begin{eqnarray*}
\bq^{\mathrm{mnl}}(\bmu) = \arg\max \displaystyle \left\{ \bmu^T\bx
- \eta \sum_{i=1}^n x_i \log x_i \ \Big{|} \ \bx \in \Delta_{n-1}
\right\}.
\end{eqnarray*}
\cite{hofbauer} further extend the result to general random utility
models. They show that for any random utility model with
continuously distributed random utility, there exists a representative agent model
that gives the same choice probability. The precise statement of
their result is as follows:

\begin{proposition}\label{prop:hof}
Let $ \bq(\bmu): {\mathcal R}^{n} \mapsto \Delta_{n-1}$ be the
choice probability function defined in (\ref{rum_choice}) where the
random vector $\bepsilon$ admits a strictly positive density on
${\mathcal R}^n$ and the function $\bq(\bmu)$ is continuously
differentiable. Then there exists $V(\cdot)$ such that:
\begin{eqnarray*}
\bq(\bmu)   = \arg\max \left\{\bmu^T\bx - V(\bx) \ \Big{|} \ \bx \in
\Delta_{n-1} \right\}.
\end{eqnarray*}
\end{proposition}

They also show that the reverse statement of Proposition
\ref{prop:hof} is not true:
\begin{proposition} [Proposition 2.2 in \citealt{hofbauer}] \label{prop:hof2} When $n
\ge 4$, there does not exist a random utility model that is
equivalent to the representative agent model with $V(\bx) = -
\sum_{i=1}^n \log{x_i}$.
\end{proposition}

Based on the above two propositions, we know that the representative
agent model strictly subsumes the random utility model as a special
case.

\subsection{Semi-Parametric Choice Model}
\label{subsec:semiparametric}

Recently, a new class of semi-parametric choice model (SCM) was
proposed by \cite{natarajan2009persistency}. Unlike the random
utility model where a certain distribution of the random utility
$\bepsilon$ is specified, in the semi-parametric choice model, one
considers a set of distributions $\Theta$ for $\bepsilon$. Given the
deterministic utility vector $\bmu$, one defines the maximum
expected utility function $w^s(\bmu)$ as follows:
\begin{eqnarray}\label{semi}
w^s(\bmu)= \displaystyle  \sup_{\theta \in \Theta} \
\mathbb{E}_{\bepsilon\sim\theta} \left[\max_{i \in \mathcal{N}} \ \
\mu_i + \epsilon_i  \right].
\end{eqnarray}
Note that in the random utility model, the maximum expected utility
function can be defined in a similar way, but only with a single
distribution $\theta$. Thus the semi-parametric choice model can be
viewed as an extension of the random utility model. Let
$\theta^*(\bmu)$ denote the extreme distribution (or a limit of a
sequence of distributions) that attains the optimal solution in
(\ref{semi}). The choice probability for alternative $i$ under this
model is given by (provided it is well-defined):
\begin{eqnarray}\label{semichoice}
q^s_i(\bmu)=  \mathbb{P}_{\theta^*(\bmu)} \left(i = \underset{k\in
\mathcal{N}}{\textrm{argmax}} \  (\mu_k + {\epsilon}_k) \right).
\end{eqnarray}

Several special cases of semi-parametric choice models have been
studied recently. One such model, called the marginal distribution
model (MDM), is proposed by \cite{natarajan2009persistency}.
In the MDM, the distribution set $\Theta$
contains all the distributions that have certain marginal
distributions. The following proposition proved in
\cite{natarajan2009persistency} shows that the marginal distribution
model can be equivalently represented by a representative agent
model:
\begin{proposition}
\label{prop:MDM} Suppose $\Theta = \left\{\theta | \epsilon_i \sim
F_i(\cdot), \forall i\right\}$ where $F_i(\cdot)$s are given
continuous distributions. Then we have:
\begin{eqnarray}\label{MDMrepre}
w^s(\bmu)= \displaystyle \max_{\bx}&\left\{\left.{\displaystyle \bmu^T\bx
+\sum_{i=1}^n \int_{1-x_i}^{1} F_i^{-1}(t)dt }\right|\bx \in
\Delta_{n-1} \right\}.
\end{eqnarray}
Furthermore, the choice probabilities $\bq^s(\bmu)$ can be obtained
as the optimal solution $\bx^*$ in (\ref{MDMrepre}).
\end{proposition}

Another semi-parametric model is the marginal moment model (MMM), in
which only the first and second moments of the marginal
distributions are known and $\Theta$ comprises all  distributions
that are consistent with the marginal moments.
\cite{natarajan2009persistency} show that the MMM
can also be represented as a representative agent model (without
loss of generality, we assume that the marginal mean of $\epsilon_i$
is $0$ for all $i$):
\begin{proposition}
\label{prop:MMM} Suppose the marginal variance of $\epsilon_i$ is
$\sigma_i$ for all $i$. Then we have
\begin{eqnarray}\label{MMMrepre}
w^s(\bmu)=\displaystyle \max_{\bx}&\left\{\left.{\displaystyle \bmu^T\bx
+\sum_{i =1}^n \sigma_i \sqrt{x_i(1-x_i)}}\right| \bx
\in\Delta_{n-1}\right\}.
\end{eqnarray}
Furthermore, the choice probabilities $\bq^{s}(\bmu)$ can be
obtained as the optimal solution $\bx^*$ in (\ref{MMMrepre}).
\end{proposition}

In order to incorporate covariance information, \cite{Mishra2010}
further propose a complete moment model (CMM), in which $\Theta$ is
the set of distributions with known first and second moments
$\bSigma$ (covariance matrix). It is shown in
\cite{ahipasaoglu2013convex} that the CMM model can also be written
as a representative agent model (again without loss of generality,
we assume the first moments are $0$):
\begin{proposition}
\label{prop:CMM} Assume $\bSigma\succ 0$. Then we have:
\begin{eqnarray}\label{CMMrepre}
w^s(\bmu) = \max_{\bx} \left\{\left.\bmu^T \bx +
\mbox{trace}\left(\bSigma^{1/2} S({\bx}) \bSigma^{1/2} \right)^{1/2}
\right| \bx \in \Delta_{n-1} \right\},
\end{eqnarray}
where $S(\bx)= \mbox{Diag}(\bx)- \bx \bx^T$ and $\mbox{trace}(X)$ is
the trace of the matrix $X$. Furthermore, the choice probabilities
$\bq^s(\bmu)$ can be obtained as the optimal solution $\bx^*$ in
(\ref{CMMrepre}).
\end{proposition}

Thus, all semi-parametric models  studied so far can be
represented as representative agent models. In the next section, we
will show that this is generally the case. Moreover, we show that
in fact, the set of representative agent models is equivalent to
that of semi-parametric  models.

Before we end this section, we comment that there are other types of
choice models in the literature in addition to those mentioned
above, such as the Markov chain-based choice model (see
\citealt{Blanchet}), the two-stage choice model (see
\citealt{Jagabathula}), the generalized attraction model (see
\citealt{gallego2014general}) and the non-parametric model (see
\citealt{farias2013nonparametric}). However, they are based on
different ideas and are less related to our study. Therefore, we
choose not to include a detailed review of those models in this
paper.

\section{Welfare-Based Choice Model}
\label{sec:unified_framework}

In this section, we propose a new framework for discrete choice
models and show that it provides a way  to unify the various choice
models reviewed in Section \ref{sec:review}. To introduce our new
model, we first notice that although various choice models
reviewed in Section \ref{sec:review} are based on different ideas,
they are all essentially functions from a vector of utilities
$\bmu$ to a vector of choice probabilities $\bq(\bmu)$. Moreover,
each of these models allows a welfare function $w(\bmu)$ that
captures the expected utility that an individual can get from the
choice model, and the choice probability vector can  be viewed as
the gradient of $w(\bmu)$ with respect to $\bmu$. Our proposed model
is based on these observations. In particular, we first construct a
class of welfare functions by defining what properties such
functions should satisfy. Then we discuss the relation between
our model and the previous ones. We start by making the following
definition:
\begin{definition}[Choice Welfare
Function]\label{def:social_welfare} Let $w(\bmu)$ be a mapping from
${\mathcal R}^n$ to $\bar{\mathcal R}$. We call $w(\bmu)$ a
choice welfare function if $w(\bmu)$ satisfies the following properties:
\begin{enumerate}
\item (Monotonicity): For any $\bmu_1$, $\bmu_2\in {\mathcal R}^n$ and $\bmu_1 \ge \bmu_2$, $w(\bmu_1) \ge
w(\bmu_2)$;
\item (Translation Invariance): For any $\bmu\in {\mathcal R}^n$, $t\in {\mathcal R}$, $w(\bmu + t {\mathbf e}) =
w(\bmu) + t$;
\item (Convexity): For any $\bmu_1$, $\bmu_2 \in {\mathcal R}^n$ and $0\le \lambda \le
1$, $\lambda w(\bmu_1) + (1-\lambda) w(\bmu_2) \ge w(\lambda\bmu_1 +
(1-\lambda) \bmu_2)$.
\end{enumerate}
In addition to the three properties, if $w(\bmu)$ is also
differentiable, then we call $w(\bmu)$ a differentiable
choice welfare function.
\end{definition}

Here we make a few comments on the three conditions in Definition
\ref{def:social_welfare}. The monotonicity condition is
straightforward. It requires that the welfare is higher if all
alternatives have higher deterministic utilities. The translation
invariance property requires that if the deterministic utilities of
all alternatives increase by a certain amount $t$, then the choice
welfare function will increase by the same amount. This is
reasonable given that choice is about relative preferences,
therefore, increasing the utilities of all alternatives by the same
amount will not change the relative preferences but will only
increase the welfare by the amount of the increment. Later, we
will see that this condition is necessary to guarantee well-defined
choice probabilities. The last condition of convexity basically
states that the welfare is higher when there is an alternative with
high utility rather than several mediocre alternatives. This is
again plausible in reality and as we will see later, all
previously reviewed choice models satisfy this condition.

In the following, we show that a choice welfare function has two
equivalent representations: a convex optimization representation and
a semi-parametric representation. This result will be instrumental
for us to derive the relations among choice models.
\begin{theorem}\label{thm:utility_equivalence} The following statements are equivalent:
\begin{enumerate}
\item [1] $w(\bmu)$ is a choice welfare function;
\item [2] There exists a convex function $V(\bx): \Delta_{n-1} \mapsto \bar{\mathcal R}$ such that
\begin{eqnarray}\label{convexdefinition}
w(\bmu) = \max \left\{\left.\bmu^T\bx - V(\bx)\right| \bx \in
\Delta_{n-1}\right\};
\end{eqnarray}
\item [3] There exists a distribution set $\Theta$ such that
\begin{eqnarray}\label{semiparametricdefinition}
w(\bmu) = \sup_{\theta\in \Theta} \mbox{
 } {\mathbb E}_{\bepsilon\sim\theta} \left[\max_{i \in {\mathcal N}} \mbox{
 } \ \mu_i + \epsilon_i\right].
\end{eqnarray}
\end{enumerate}
\end{theorem}

{\noindent\bf Proof.} First we show that the $w(\bmu)$ defined in
(\ref{convexdefinition}) and (\ref{semiparametricdefinition}) are
choice welfare functions. To see this, we note that the monotonicity
and translation invariance properties are immediate from
(\ref{convexdefinition}) and (\ref{semiparametricdefinition}). For
the convexity, we note that $w(\bmu)$ defined in
(\ref{convexdefinition}) is the supremum of linear functions of
$\bmu$ thus is convex in $\bmu$. In
(\ref{semiparametricdefinition}), for each $\bepsilon$, $\max_{i\in
\mathcal{N}} \left\{\mu_i + \epsilon_i\right\}$ is a convex function
in $\bmu$, and so is the expectation. Therefore, if $w(\bmu)$ is defined
by (\ref{convexdefinition}) or (\ref{semiparametricdefinition}),
then it must be a choice welfare function.

Next we show the other direction. That is, if $w(\bmu)$ is a
choice welfare function, then it can be represented in the form of
(\ref{convexdefinition}) and (\ref{semiparametricdefinition}). First
we note that if a choice welfare function $w(\bmu) = \infty$ for
some $\bmu$, then by the translation invariance property and the
monotonicity property, it must be that $w(\bmu) = \infty$ for any
$\bmu$. In that case, we can choose $V(\bx) = -\infty$ and $\Theta =
\{\theta_\infty\}$ where $\theta_\infty$ is a singleton distribution
taking value on $(\infty,...,\infty)$. Therefore, $w(\bmu)$ can be
represented by (\ref{convexdefinition}) and
(\ref{semiparametricdefinition}) in that case. Similarly, if
$w(\bmu) = -\infty$ for some $\bmu$, then it must be that $w(\bmu) =
-\infty$ for all $\bmu$, and we can take $V(\bx) = \infty$ and
$\Theta = \{\theta_{-\infty}\}$, where $\theta_{-\infty}$ is a
singleton distribution on $(-\infty,...,-\infty)$. Therefore,
$w(\bmu)$ can be represented in (\ref{convexdefinition}) and
(\ref{semiparametricdefinition}) in this case too.

In the remainder of the proof, we focus on the case where $w(\bmu)$
is finite for all $\bmu$. In this case, by Proposition 1.4.6 of
\cite{bertsekas}, $w(\bmu)$ must be continuous. The remaining proof
is divided into two parts:

1. We show that any choice welfare function $w(\bmu)$ can be
represented by (\ref{convexdefinition}). Since $w(\bmu)$ is monotone
and translation invariant, the following holds:
\begin{eqnarray*}
w(\bmu) = \min_{\by} \left\{\ w(\by) + \max_{i} \left\{\mu_i -
y_i\right\}\right\} =  \min_{\by}\ \left\{ w(\by) + \max_{\bx \in
\Delta_{n-1}} (\bmu - \by)^T \bx\right\}.
\end{eqnarray*}
Here the first equality holds since for any $\by$, $w(\bmu) = w(\bmu -
\max_i\left\{\mu_i - y_i\right\}\mathbf{e}) +  \max_{i} \left\{\mu_i
- y_i\right\}$ by the
translation invariance property.
Furthermore, by the monotonicity property, $w(\bmu -
\max_i\left\{\mu_i - y_i\right\}\mathbf{e}) \le w(\by)$ and
the equality holds when $\by = \bmu$.

Next we define $L(\bx, \by) = w(\by) + (\bmu - \by)^T \bx$. We have
for fixed $\bx$, $L(\bx, \cdot)$ is convex in $\by$ (by the
convexity of $w(\cdot)$); and for fixed $\by$, $L(\cdot, \by)$ is
convex and closed in $\bx$. Furthermore, $\inf_{\by} \ \max_{\bx\in
\Delta_{n-1}} L(\bx, \by) = w(\bmu) < \infty$ and the function
$p(\bu) = \inf_{\by} \ \max_{\bx\in \Delta_{n-1}} \left\{L(\bx, \by) -
\bu^T \bx\right\} = w(\bmu - \bu)$ is continuous. Therefore, by
Proposition 2.6.2 of \cite{bertsekas}, the minimax equality holds,
i.e.,
\begin{eqnarray*}
\inf_{\by} \max_{\bx\in \Delta_{n-1}} L(\bx, \by) = \max_{\bx\in
\Delta_{n-1}}\inf_{\by} L(\bx, \by).
\end{eqnarray*}
Therefore, we have:
\begin{equation*}
 w(\bmu)=\max_{\bx \in \Delta_{n-1}} \left\{\bmu^T\bx + \inf_{\by} \left\{w(\by)  - \by^T\bx \right\}\right\}=\max_{\bx \in \Delta_{n-1}} \{ \bmu^T \bx -
 V(\bx)\}
\end{equation*}
where $V(\bx) =  \sup_{\by }\{ \by^T \bx - w(\by)\}$ is a convex
function.

2. Next we show that any choice welfare function can be represented
by (\ref{semiparametricdefinition}). Since $w(\bmu)$ is convex,
there exists a subgradient for any $\bmu$. We denote the subgradient
vector by $\bd(\bmu) = (d_1(\bmu),\dots, d_n(\bmu))^{T}$. Here it is
possible that the choice of $\bd(\bmu)$ is not unique, in that case,
we can choose an arbitrary one. Furthermore, by taking the
derivative with respect to $t$ in the translation invariance
equation, and by applying the chain rule (see Proposition 4.2.5 of
\citealt{bertsekas}), we have for any subgradient $\bd(\bmu)$, it
must hold that ${\mathbf e}^T\bd(\bmu) = 1$. Similarly, by the
monotonicity property of $w(\bmu)$, we must have $\bd(\bmu)\ge
\mathbf{0}$. By the definition of subgradient and the convexity of
$w(\bmu)$, we must have:
\begin{eqnarray*}
w(\bmu) \ge (\bmu - \bz)^T \bd(\bz) + w(\bz), \quad\forall \bz \in \mathcal{R}^{n},
\end{eqnarray*}
where the equality holds when $\bz = \bmu$. Define $l(\bz) = w(\bz)
- \bz^T\bd(\bz)$. By reorganizing terms, we have
\begin{eqnarray}\label{proofstep1}
w(\bmu) = \sup_{\bz} \{\bmu^T \bd(\bz) + l(\bz)\}.
\end{eqnarray}

Now we define the distribution set as follows: Let $\Theta =
\{\theta_{\bz} \big{|} \bz\in{\mathcal R}^n\}$, where $\theta_{\bz}$
is an $n$-point distribution with
\begin{eqnarray*}
{\mathbb P}_{\theta_{\bz}}\left(\bepsilon = \bepsilon_{\bz}^i\right)
= d_i(\bz), \quad\mbox{for }i = 1,...,n
\end{eqnarray*}
where
\begin{eqnarray*}
\epsilon_{\bz}^i(j) = \left\{\begin{array}{ll} l(\bz) & \mbox{if }
j = i\\
-\infty & \mbox{if } j\neq i.
\end{array}
\right.
\end{eqnarray*}
That is, $\epsilon_{\bz}^i$ is a vector of all $-\infty$'s except
$l(\bz)$ at the $i$th entry. Therefore, for any $\bz$, we have
\begin{eqnarray*}
\mathbb{E}_{\theta_{\bz}} [\max_{i} \ \left\{ \mu_i +
\epsilon_i\right\}]  =  \sum_{i=1}^n d_i(\bz) (\mu_i + l(\bz)) =
\bmu^T \bd(\bz) + l(\bz).
\end{eqnarray*}
Then by (\ref{proofstep1}), we have
\begin{eqnarray*}
w(\bmu) = \sup_{\bz} \mbox{ }\{\bmu^T\bd(\bz) + l(\bz)\} =
\sup_{\bz} \mbox{ }{\mathbb E}_{\theta_{\bz}}[\max_i \ \left\{\mu_i
+ \epsilon_i\right\}] = \sup_{\theta\in \Theta}\mbox{ } {\mathbb
E}_{\theta} [\max_i \ \left\{\mu_i + \epsilon_i\right\}].
\end{eqnarray*}
Therefore, the theorem is proved. $\hfill\Box$

Note that the above discussion focuses on the equivalent representations of the choice welfare
function. In the following we establish its implication to discrete choice models. In this paper, we refer to discrete choice models as the entire set of functions $\bq(\bmu): \mathcal{R}^{n} \mapsto \Delta_{n-1}$, mapping a utility vector to a choice probability vector. We first
propose the following choice model based on the choice welfare
function:
\begin{definition}[Welfare-based Choice Model] Suppose $w(\bmu)$ is
a differentiable choice welfare function. Then the welfare-based
choice model derived from $w(\bmu)$ is defined by
\begin{eqnarray}\label{socialchoice}
\bq(\bmu) = \nabla w(\bmu).
\end{eqnarray}
\end{definition}

Note that when $w(\cdot)$ is differentiable, we have $\nabla w(\bmu)
\in \Delta_{n-1}$ by the translation invariance property of
$w(\bmu)$. Therefore $\bq(\bmu)$ defined by (\ref{socialchoice}) is
indeed a valid choice model. Next we show the equivalence of various
choice models. We first introduce the following definitions (see
\citealt{rockafellar}):

\begin{definition}[Proper Function] A function $f : X \mapsto \bar{\mathcal R}$ is proper
if $f(\bx) < \infty$ for at least one $\bx \in X $ and $f(\bx) >
-\infty$ for all $\bx\in X$.
\end{definition}

\begin{definition}[Essentially Strictly Convex Function]
A proper convex function $f$ on ${\mathcal R}^n$ is essentially
strictly convex if $f$ is strictly convex on every convex subset of
\begin{eqnarray*}
\mbox{dom} (\partial f) = \left\{\bx\big{|} \partial f(\bx) \neq
\phi\right\}.
\end{eqnarray*}
where $\partial f(\bx)$ is the set of subgradients of $f$ at $\bx$, and $\phi$ is the empty set.
\end{definition}

Note that any strictly convex function is essentially strictly
convex. Next we have the following theorem, whose proof is relegated
to the Appendix:
\begin{theorem}\label{thm:choicemodel_equivalence}
For a choice model $\bq: {\mathcal R}^n \mapsto \Delta_{n-1}$, the
following statements are equivalent:
\begin{enumerate}
\item There exists a differentiable choice welfare function $w(\bmu)$ such that
$\bq(\bmu) = \nabla w(\bmu)$;
\item There exists an essentially strictly convex function $V(\bx)$ such that
\begin{eqnarray*}
\bq(\bmu) = \arg\max\left\{\bmu^T\bx - V(\bx) \Big{|}
\bx\in\Delta_{n-1}\right\};
\end{eqnarray*}
\item There exists a distribution set $\Theta$ such that
\begin{eqnarray*}
\bq(\bmu) = \nabla_{\bmu} \left\{\sup_{\theta\in \Theta} \mbox{
 } {\mathbb E}_\theta \left[\max_{i\in{\mathcal N}} \mbox{
 }\mu_i + \epsilon_i\right]\right\}.
\end{eqnarray*}
\end{enumerate}
\end{theorem}

The next corollary follows immediately from Theorem
\ref{thm:choicemodel_equivalence} and Propositions \ref{prop:hof}
and \ref{prop:hof2}.
\begin{corollary}\label{cor:rum}
Let $\bq(\bmu)$ be a random utility model with absolutely continuous
distribution $\theta$ and $w(\bmu)$ be the corresponding expected
utility an individual can get under this model. Then $w(\bmu)$ is a
differentiable choice welfare function, and $\bq(\bmu) = \nabla
w(\bmu)$. Moreover, the reverse statement is not true, i.e., there
exists a differentiable choice welfare function $w(\bmu)$ such that
there is no random utility model that gives the choice probability
$\bq(\bmu) = \nabla w(\bmu)$.
\end{corollary}

The significance of Theorems \ref{thm:utility_equivalence} and
\ref{thm:choicemodel_equivalence} is mainly twofold. First, we
propose a new framework for discrete choice model, the welfare-based
choice model, which is based on the desired functional properties of
the expected utility function. With the help of the
new framework, we establish the connection between two
existing choice models, the representative agent model and the
semi-parametric model. In particular, we show that those two classes
of choice models are equivalent. This result explains the prior
results that for every known semi-parametric model, there is a
corresponding representative agent model. In addition, it asserts
that the reverse is also true, which is quite surprising in some
sense. Therefore, in terms of the scope of choice models that can be
captured, those three models (the welfare-based choice model, the
representative agent model and the semi-parametric model) are the
same. We believe this result is useful for the theoretical study of
discrete choice models.

Second, by establishing the equivalence of the three classes of
choice models, we can allow more versatile ways to construct a
choice model. In particular, we can pick any of the three
representations to start with. For the welfare-based choice model,
one needs to choose  a choice welfare function $w(\bmu)$
which satisfies the three conditions. For the representative agent
model,  one needs to choose a (strictly) convex
regularization function. And for the semi-parametric model,  one
needs to choose  a set of distributions. In different situations,
it might be easier to use one representation than the other in order
to capture certain properties of the choice model. In addition, by
Corollary \ref{cor:rum}, the welfare-based choice model strictly
subsumes the random utility model, thus it is possible to construct
new choice models that have certain interesting properties that a
random utility model could not accommodate. We will further study
this issue in Sections \ref{sec:relation_rum} and
\ref{sec:substitutability}.

The next theorem studies one desirable
property of choice models and investigates how it can be reflected to the
construction of the three choice models. We start with the following
definition:
\begin{definition}[superlinear choice welfare function]
\label{def:boundedwelfare} A differentiable choice welfare function
$w(\bmu)$ is called \textit{superlinear} if there
exist $b_i, i=1,...,n$, such that for any $\bmu \in {\mathcal
R}^n$:
\begin{eqnarray*}
w(\bmu) \geq \mu_i + b_i, \ \forall \  i = 1,...,n.
\end{eqnarray*}
\end{definition}

This property is desirable in most applications. It requires that
the utility one can get from a set of alternatives is not much less
than the utility of each alternative. After all, for each
alternative $i$, one can always choose it and obtain the
corresponding utility. We have the following theorem:
\begin{theorem}\label{thm:boundequivalence}
For a choice model $\bq: {\mathcal R}^n \mapsto \Delta_{n-1}$, the
following statements are equivalent:
\begin{enumerate}
\item There exists a superlinear differentiable choice welfare function $w(\bmu)$ such that
$\bq(\bmu) = \nabla w(\bmu)$;
\item There exists an essentially strictly convex function $V(\bx)$ that is upper bounded on $\Delta_{n-1}$ such that
\begin{eqnarray*}
\bq(\bmu) = \arg\max\left\{\bmu^T\bx - V(\bx) \Big{|}
\bx\in\Delta_{n-1}\right\};
\end{eqnarray*}
\item There exists a distribution set $\Theta$ containing only
distributions with finite expectation (i.e., ${\mathbb E}_{\theta} |\epsilon_i| < \infty$ for all
$i$ and $\theta \in \Theta$) such that
\begin{eqnarray*}
\bq(\bmu) = \nabla_{\bmu} \left\{\sup_{\theta\in \Theta} \mbox{
 } {\mathbb E}_\theta \left[\max_{i\in{\mathcal N}} \mbox{
 }\mu_i + \epsilon_i\right]\right\},
\end{eqnarray*}

\end{enumerate}
Moreover, if either of the above cases holds, then $\bq(\bmu)$ can
span the whole simplex, i.e., for all $\bx$ in the interior
of $\Delta_{n-1}$, there exists $\bmu$ such that $\bq(\bmu) = \bx$.
\end{theorem}

We present the proof of Theorem \ref{thm:boundequivalence} in the
Appendix. We can see that Theorem \ref{thm:boundequivalence} further
develops the equivalence of choice models obtained in Theorem
\ref{thm:choicemodel_equivalence} by narrowing down the discussion
to welfare-based choice models with the desirable superlinear
property. In particular, we find that a superlinear differentiable
choice welfare function has a semi-parametric representation, of
which the distribution set contains at least one bounded
distribution. The distribution set containing bounded distribution
is also desirable due to its potential practical application. The
last statement that $\bq(\bmu)$ spans the whole simplex is related
to the results in \cite{hofbauer}, \cite{norets2013surjectivity} and
\cite{mishra2014theoretical}. These papers provide  conditions
under which $\bq(\bmu)$ defined from the RUM or the MDM can span the whole
simplex. Theorem \ref{thm:boundequivalence} extends these results to
more general conditions.

\section{Relation to the Random Utility Model} \label{sec:relation_rum}

In the last section, we proposed a new framework for choice models:
the welfare-based choice model. In particular, by Corollary
\ref{cor:rum}, the class of welfare-based choice models strictly subsumes the
random utility model. In this section, we investigate further the
relation between the welfare-based choice model and the random
utility model. In particular, we study under what conditions a
welfare-based choice model can be equivalently represented by a
random utility model. This study will help us understand clearly the
relations between various choice models and the random utility model
and design new choice models that do not necessarily have a random
utility representation.

First, we show that when there are only two alternatives, the class
of random utility models is equivalent to the class of welfare-based
choice models.

\begin{theorem}\label{thm:RUMequivalence}
For any differentiable choice welfare function $w(\mu_1,\mu_2)$,
there exists a distribution $\theta$ of $\{\epsilon_1, \epsilon_2\}$
such that:
\begin{eqnarray}\label{rumequivalence}
w(\mu_1,\mu_2) = \mathbb{E}_{\theta} [\max\{\mu_1 + \epsilon_1,
\mu_2 + \epsilon_2\}].
\end{eqnarray}
In addition, if $w(\mu_1,\mu_2)$ is superlinear, then there exists a
 distribution $\theta$ with finite expectation (i.e., ${\mathbb E}_{\theta} |\epsilon_1| < \infty$ and ${\mathbb E}_{\theta} |\epsilon_2| < \infty$) that satisfies (\ref{rumequivalence}).
\end{theorem}
\noindent{\bf Proof.} Define  $v(x)\triangleq w(x,0)$. Since
$w(\cdot)$ is differentiable,  by the chain rule, we have
$$v'(x)  = {\partial w  \over \partial \mu_1 } (x,0).$$
Since $w(\mu_1,\mu_2)$ is convex and satisfies the translation
invariance property, we have $v'(x) \in [0,1] $ and is increasing.
We define a distribution $\theta$ of $\{\epsilon_1, \epsilon_2\}$ as
follows:
$$\{\epsilon_1, \epsilon_2\} = \big\{v_0 - \max\{\xi,0\},v_0- \max\{-\xi,0\}\big\}, $$
where $v_0= v(0)=w(0,0)$ and $\xi$ is a random variable with c.d.f.
$F_{\xi}(x) = {\mathbb P}(\xi \le x) = v'(x)$. Note $F(\cdot)$ is a
well-defined c.d.f. since $w(\cdot)$ is convex and differentiable,
thus $v'(x)$ must be continuous and increasing
(\citealt{rockafellar}).

Now we compute $\mathbb{E}_{\theta} [\max\{\mu_1 + \epsilon_1, \mu_2
+ \epsilon_2\}]$. We have
\begin{eqnarray*}
{\mathbb E}_{\theta} [\max\{\mu_1 + \epsilon_1, \mu_2 +
\epsilon_2\}]  & = &  \mu_1 + v_0 + {\mathbb E}_{\theta}[
\max\{-\max\{\xi, 0\}, \mu_2 - \mu_1 - \max\{-\xi, 0 \}\}] \\
& = & \mu_1 + v_0 + {\mathbb E}_{\theta}[ \max\{0, \mu_2 - \mu_1 +
\xi\} -\max\{\xi,0\}],
\end{eqnarray*}
where the last step can be verified by considering $\xi \ge 0$ and
$\xi \le 0$, respectively.

Now we compute the last term. For $x\ge 0$, we have (let $\mathbb
I(\cdot)$ be the indicator function):
\begin{eqnarray*}
{\mathbb E}_{\theta} [ \max\{0, x + \xi\} - \max\{0,  \xi\}]
& = & x {\mathbb P} ( \xi > 0) + {\mathbb E}_{\theta} [(x + \xi) \cdot {\mathbb I}(-x < \xi \le 0) ] \\
& = & x {\mathbb P} ( \xi > 0) + \int_{-x}^{0} (x+ \xi) dv'(\xi) \\
& = & x (1- v'(0)) + (x + \xi) v'(\xi)\mid_{-x}^0 - \int_{-x}^0 v'(\xi) d\xi \\
& = & x - v_0 + v(-x).
\end{eqnarray*}
Similarly, for $x\le 0$, we have
\begin{eqnarray*}
{\mathbb E}_{\theta} [ \max\{0, x + \xi\} - \max\{0,  \xi\}]
& = & x {\mathbb P} ( \xi > -x) + {\mathbb E}_{\theta} [-\xi \cdot {\mathbb I}(0 < \xi \le -x) ] \\
& = & x {\mathbb P} ( \xi > -x ) - \int_{0}^{-x} \xi dv'(\xi) \\
& = & x (1- v'(-x))-  \xi v'(\xi) \mid_{0}^{-x} + \int_{0}^{-x}v'(\xi)d\xi \\
& = & x - v_0 + v(-x).
\end{eqnarray*}
Therefore, $\mathbb{E}_{\theta} [\max\{\mu_1 + \epsilon_1, \mu_2 +
\epsilon_2\}]= \mu_1 + v_0 + (\mu_2 - \mu_1) - v_0 + v(\mu_1 -
\mu_2) = w(\mu_1, \mu_2)$.

To prove the last statement, it suffices to show that both ${\mathbb
E}_{\theta} [ \max\{0, \xi\}]$ and ${\mathbb E}_{\theta} [ \max\{0,
-\xi\}]$ are finite if $w(\bmu)$ is superlinear. If $w(\cdot)$ is
superlinear, then we have $v(t)-t = w(0,-t) $ is decreasing in $t$
and lower bounded, thus $L_1 = \lim_{t \rightarrow + \infty}
(v(t)-t) $ exists and is finite. Similarly, $v(t) = w(t,0) $ is
increasing in $t$ and lower bounded, thus $L_2 = \lim_{t \rightarrow
- \infty} v(t)$ exists and is finite. Therefore, we have:
\begin{eqnarray*}
{\mathbb E}_{\theta} [ \max\{0, \xi\}] = \int_{0}^{+\infty} {\mathbb
P}_{\theta} \left(\xi \geq t\right) dt= \int_{0}^{+\infty}( 1-v'(t))
dt = (t-v(t))\big|^{+\infty}_{0}
   = v(0) - L_1,
\end{eqnarray*}
and
\begin{eqnarray*}
{\mathbb E}_{\theta} [ \max\{0, -\xi\}] = \int_{0}^{+\infty}
{\mathbb P}_{\theta} \left(-\xi \geq t\right) dt =
\int_{0}^{+\infty}  v'(-t)  dt = \int_{-\infty}^{0}  v'(t)  dt =
v(0) - L_2.
\end{eqnarray*}
Thus, the theorem is proved. $\hfill\Box$

%So far we have established the equivalence of these three classes of choice models: the welfare-based choice model, the representative agent choice model and the semi-parametric choice model. The equivalence among the three classes of choice models implies that we can begin with any of the three representations to construct the same choice model. In the remainder of this section, we provide a few examples of choice models that have an explicit welfare-based representation. Among them, there are several new choice models that might be of interest.
%\subsection{Relationship with the Generalized Extreme Value (GEV)
% model}

By Proposition \ref{prop:hof2}, when $n \ge 4$, the welfare-based
choice model strictly subsumes the random utility model. In fact, as
we will see in some examples later (Examples \ref{example:2} and
\ref{example:1} in Section \ref{sec:substitutability}), this is also
true for $n = 3$. In light of this relation between these two
classes of choice models, it would be interesting to know the exact
difference between them. In other words, it would be interesting to
know what property is restricted in the random utility model but not
in the welfare-based choice model. In the following, we pinpoint
this difference. The following result is a direct consequence of the
result in \citet{McFadden_nested}:
\begin{proposition}\label{propw_rum}
Let $w(\bmu):{\mathcal R^n} \mapsto {\mathcal R}$ be a
differentiable function. Then $\nabla w(\bmu)$ is consistent with a
random utility model if and only if $w(\cdot)$ satisfies the
monotonicity, translation invariance, convexity properties, and for
any $k \ge 1$ and $i_1,...,i_k$ all distinct,
\begin{eqnarray*}
(-1)^k \frac{\partial^k w(\bmu)}{\partial
\mu_{i_1},...,\partial\mu_{i_k}} \le 0.
\end{eqnarray*}
\end{proposition}

By Proposition \ref{propw_rum} and the above discussions, we point
out that the difference between a random utility model and a
welfare-based choice model (thus also the representative agent model
and the semi-parametric model by Theorem
\ref{thm:choicemodel_equivalence}) lies in the requirement on the
higher-order derivatives of $w(\cdot)$. In particular, a random
utility model requires that the higher-order cross-partial
derivatives of $w(\cdot)$ have alternating signs, while in the
welfare-based choice model, it only requires that the Hessian matrix
of $w(\cdot)$ be positive semidefinite, and there is no requirement
on other higher-order derivatives. This difference will enable us to
better understand
 the difference between those models and later construct
choice models with new properties.

We next consider an important subclass of the random utility model:
the generalized extreme value (GEV) model. The GEV model was first
proposed by \cite{MacfaddenGEV}. It is a special case of the random
utility model in which the random part of the utility $\epsilon_i$s
take a joint generalized extreme value distribution. The GEV model
covers various popular models, including the MNL model, the nested
logit model, etc. An equivalent definition of the GEV model is given
as follows (\citealt{McFadden_nested}):
\begin{definition}[GEV model] \label{def:gev}
A choice model $\bq(\bmu)$ is a GEV model if and only if there
exists a function $H(\by): \mathcal{R}_{+}^{n}\mapsto \mathcal{R}$
such that
\begin{eqnarray}\label{choiceprobgev}
\bq(\bmu)=\eta \nabla_{\mbs{\mu}}\mathrm{log}
H(e^{\mu_{1}},\dots,e^{\mu_{n}}),
\end{eqnarray}
where $H(\by)$ satisfies the following properties:
\begin{enumerate}
\item $H(\mb{y}) \geq 0$ for all $\mb{y} \in \mathcal{R}_{+}^{n}$.
\item $H(\mb{y})$ is homogeneous of degree $1/\eta$, i.e., $H(\alpha\mb{y}) = \alpha^{1/\eta} H(\mb{y})$.
\item $H(\mb{y}) \rightarrow \infty$ as $y_j \rightarrow \infty$ for any $j$.
\item The $k$th-order cross-partial derivatives of $H(\mb{y})$ exist for
all $1\le k\le n$, and for all distinct $i_1,...,i_k$,
\begin{eqnarray*}
(-1)^k { \partial^k H(\mb{y}) \over \partial y_{i_1} ...\partial
y_{i_k}} \leq 0.
\end{eqnarray*}
\end{enumerate}
\end{definition}
%for which
%% Then, $F(z_{1},\dots,z_{n})=\exp[-H(e^{-z_{1}},\dots,e^{-z_{n}})]$ is a multivariable distribution function, and the choice probabilities that result from the maximization of the random utilities for which the multivariate distribution function is given by $F(\cdot)$ are equal to
%A choice model is called a GEV model if and only if there exists
%such a representation.
%\begin{equation*}
%  q_i(\bmu)= { {e^{\mu_{i}} H^{(1)}_i(e^{\mu_{1}},\dots,e^{\mu_{n}})} \over H(e^{\mu_{1}},\dots,e^{\mu_{n}})}, \quad \forall i \in \mathcal{N}.
%\end{equation*}

Under appropriate specifications of $H(\cdot)$, various known choice
models can be obtained from the GEV model. We list the MNL model and
the nested logit model as examples (\citealt{train2009discrete}).
\begin{itemize}
\item MNL model. If one chooses $H(\mb{y}) = \sum_{i\in \mathcal{ N}}
y_i^{1/\eta}$, then the corresponding choice model is the MNL model
with choice probabilities:
\begin{equation*}
q_i(\bmu)= \frac{\exp({\mu_{i}/\eta})}{\underset{k\in
\mathcal{N}}\sum \exp({\mu_{k}/\eta})}.
\end{equation*}
\item Nested Logit model. Suppose the $n$ alternatives are partitioned into $K$ nests labeled $B_1,...,B_K$.
If one chooses $ H(\mb{y}) = \sum_{l=1}^K \left( \sum_{i \in B_l}
y_i^{1/\lambda_l}\right)^{\lambda_l}$, then the corresponding choice
model is the nested logit model with choice probabilities:
\begin{equation*}
q_i(\bmu)=  {\exp(\mu_i / \lambda_k) (\sum_{j \in B_k} \exp(\mu_j /
\lambda_k))^{\lambda_k -1} \over \sum_{l=1}^K  \left(\sum_{j \in
B_l} \exp(\mu_j / \lambda_l) \right)^{\lambda_l}}.
\end{equation*}
\end{itemize}

Since the welfare-based choice model strictly subsumes the random
utility model, we know that the GEV model can be equivalently
represented by welfare-based choice models. By
$\bq(\bmu)=\nabla_{\bmu}w(\bmu)=\eta\nabla_{\bmu}\log
H(e^{\mu_{1}},\dots,e^{\mu_{n}})$, it is implied that a GEV model
derived from a specific $H(\cdot)$ is equivalent to a welfare-based
choice model with welfare function $w(\bmu)=\eta \mathrm{log}
H(e^{\mu_{1}},\dots,e^{\mu_{n}})$ and such $w(\cdot)$ must satisfy
the properties in Proposition \ref{propw_rum}.

Conversely, for any $w(\bmu)$ being a differentiable choice welfare
function, we can define a function $H(\by)=\exp(w(\log
y_{1},\dots,\log y_{n})/\eta)$. The following discussions point out
what properties such an $H(\cdot)$ would satisfy:
\begin{enumerate}
\item By definition, $H(\bz)\ge 0$ for all $\bz$.
\item Since $w(\cdot)$ is translation invariant, we have that
$H(\alpha\bz) = \alpha^{1/\eta} H(\bz)$, i.e., $H(\bz)$ is
homogeneous of degree $1/\eta$.
\item Since $w(\cdot)$ is monotone, we have
\begin{eqnarray*}
{\partial w(\bmu) \over\partial \mu_i} = { \eta \exp(\mu_i)
H^{(1)}_i(\bz) \over H(\bz)} \geq 0, \quad \forall i \in
\mathcal{N}.
\end{eqnarray*}
where $\bz = (e^{\mu_{1}},\dots,e^{\mu_{n}})$ and $H_i^{(1)}(\cdot)$
is the partial derivative of $H$ with respect to $i$. Therefore, all
first-order partial derivatives of $H(\bz)$ are non-negative.
\item Last, in order for $w(\cdot)$ to be convex, we need the
Hessian matrix of $w(\cdot)$ defined as follows to be positive
semidefinite:
\begin{eqnarray*}
{\partial^2 w(\bmu) \over \partial \mu_i \partial \mu_j}={ \eta
\exp(\mu_i + \mu_j) \left(H(\bz) \cdot H^{(2)}_{ij}(\bz) -
H^{(1)}_i(\bz) \cdot H^{(1)}_j(\bz)\right) \over H^2(\bz)}
\quad\mbox{ and } \quad {\partial^2 w(\bmu) \over \partial \mu_i^2}
= - \sum_{j \neq i} {\partial^2 w(\bmu) \over \partial \mu_i
\partial \mu_j}.
\end{eqnarray*}
where $\bz = (e^{\mu_{1}},\dots,e^{\mu_{n}})$, $H_i^{(1)}$ is the
partial derivative of $H(\cdot)$ with respect to $i$, and
$H_{ij}^{(2)}$ is the second-order partial derivative of $H(\cdot)$
with respect to $i$ and $j$.
\end{enumerate}

It is worth pointing out that the last condition holds if all
second-order cross-partial derivatives of $H$ are negative, but the reverse
is not necessarily true (the equivalent condition involves all the
zero-, first- and second-order derivatives of $H$). Therefore, the
GEV model requires an even stronger condition that the higher-order
derivatives of $\exp(w(\log{y_1},...,\log{y_n})/\eta)$ have
alternating signs, while in the welfare-based choice model, we only
need the first-order derivative to be positive and some condition
that is weaker than requiring all the cross second-order derivatives
 be negative.

\section{Substitutability and Complementarity of Choices}
\label{sec:substitutability}

In the previous section, we have seen that the distinction between
the welfare-based choice model and the random utility model lies in
the property of the higher-order derivatives of the choice welfare
function. In particular, the random utility model has stronger
requirements on the higher-order derivatives. In this section, we
will discuss more in depth about the practical meaning of such
properties. We introduce two concepts, which we call the {\it
substitutability} and {\it complementarity} of choices and then
discuss the practical relevance of these two concepts. We show that
if a choice model is derived from a random utility model, then the
alternatives can only exhibit substitutability. However, using our
welfare-based framework, we can design choice models that have more
flexible substitutability or complementarity patterns. We also show
how this property can be reflected through the regularization
function in a representative agent model. Before we formally define
these two concepts, we first introduce the definition of local
monotonicity:
\begin{definition}[local monotonicity]\label{def:localin}
A function $f(x): \mathcal{R} \mapsto \mathcal{R}$ is locally increasing at $x$ if there exists $\delta>0$ such that
$$f(x-h)\leq f(x) \leq f(x+h), \ \forall \ 0<h<\delta.$$
Similarly, $f(x)$ is locally decreasing at $x$ if there exists $\delta>0$ such that
$$f(x-h)\geq f(x) \geq f(x+h), \ \forall \ 0<h<\delta.$$
\end{definition}

Now we introduce the definition of substitutability and
complementarity in choice models:
\begin{definition}\label{def:substitutability}
Consider a choice model $\bq(\bmu): {\mathcal R}^n \mapsto
\Delta_{n-1}$. For any fixed $\bmu$ and $i, j\in {\mathcal N}$:
\begin{enumerate}
\item (Substitutability) If $q_j(\bmu)$ is locally decreasing in $\mu_i$ at $\bmu$, then we
say alternative $i$ is substitutable to alternative $j$ at $\bmu$.
Furthermore, if $q_j(\bmu)$ is locally decreasing in $\mu_i$ for all $\bmu$, then
we say alternative $i$ is substitutable to alternative $j$;
\item (Complementarity) If $q_j(\bmu)$ is locally increasing in $\mu_i$ at $\bmu$, then we
say alternative $i$ is complementary to alternative $j$ at $\bmu$.
Furthermore, if $q_j(\bmu)$ is locally increasing in $\mu_i$ for all $\bmu$, then
we say alternative $i$ is complementary to alternative $j$.
\item (Substitutable Choice Model) For all $i \neq j$, if alternative $i$ is substitutable to alternative $j$, then we say $\bq(\bmu)$ is a substitutable choice model.
\end{enumerate}
\end{definition}

The definition of substitutability and complementarity of two
alternatives is similar to that of two consumer goods (see
\citealt{mankiw}). However, in Definition
\ref{def:substitutability}, the independent variable is not the
price, and the dependent variable is choice probability rather than
demand. We first investigate some basic properties of
substitutability and complementarity.

\begin{proposition}\label{prop:basicpropertysub}
Consider a choice model $\bq(\bmu): \mathcal{R}^n \mapsto
\Delta_{n-1}$ that is derived from a differentiable choice welfare
function $w(\bmu)$. For any $i$, alternative $i$ must be
complementary to itself. Furthermore, if $w(\bmu)$ is second-order
continuously differentiable and alternative $i$ is substitutable
(complementary, resp.) to alternative $j$ at $\bmu$, then
alternative $j$ must be substitutable (complementary, resp.) to
alternative $i$ at $\bmu$.
 \end{proposition}

The proof of Proposition \ref{prop:basicpropertysub} uses some basic
properties of continuous and convex function and is delegated to the
Appendix. The proposition shows that when $w(\bmu)$ is second-order
continuously differentiable, the substitutability (complementarity,
resp.) property is a reciprocal property. In these cases, we shall
say $i$ and $j$ are substitutable (complementary, resp.) in the
following discussions.

In the following, we investigate the substitutability and complementarity of choice models. First we show that random utility models are all substitutable:
\begin{theorem}\label{thm:randomutilitysub}
Any random utility  model $\bq(\bmu)$ is a substitutable choice
model.
\end{theorem}

Theorem \ref{thm:randomutilitysub} directly follows from Proposition
\ref{propw_rum}. It states that in a random utility model, if the
utility of one alternative increases while the utilities of all
other alternatives stay the same, then it must be that the choice
probabilities of all other alternatives decrease. This is certainly
plausible in practice, especially if $\bmu$ is intepreted as how much
a consumer values each product. However, as we show in the
following example, sometimes it might be desirable to allow
different alternatives to exhibit certain degrees of
complementarity. This is especially true if we allow more versatile
meanings of the utility $\bmu$.
\begin{example}\label{example:substitutability}
Suppose a customer is considering to buy a camera from the following
three alternatives: a Canon-A model, a Canon-B model and a Sony-C
model. On a certain website, there are customer reviews that rate
each model, which we denote by $v_1$, $v_2$ and $v_3$, respectively.
We assume that the customer's choice is solely based on those review
scores (suppose other factors are fixed). That is, the choice
probability $\bq$ is a function of $\bv = (v_1,v_2,v_3)$. Suppose at
a certain time, a new review for the Canon-A model comes in, rating
it favorably. How would it change the purchase probability of the
Canon-B model?

The answer to the above questions may depend. There might be two
forces. On one hand, due to a new favorable rating  given to the
Canon brand, the probability of choosing the Canon-B model might
increase. On the other hand, the favorable rating for the Canon-A model
might switch some customers from the Canon-B model to the Canon-A model.
Either force might be dominant in practice. If the former force
is stronger, then it is plausible that one additional favorable
rating for the Canon-A model might increase the choice probability of
the Canon-B model. $\hfill\Box$
\end{example}

The above example illustrates that sometimes it might be desirable
to have a choice model in which a certain pair of alternatives exhibit
complementarity. One may notice that the above example may be
reminiscent of the nested logit model, in which the customers first
choose a nest (in this case, the brand), and then choose a
particular product. When increasing the utility of another product
in the same nest, the tradeoff is between the probability of
choosing the nest (which will be higher) and the individual product
(which will be lower). However, we note that the nested logit model
is essentially a random utility model (with the randomness
$\bepsilon$ chosen to be an extreme value distribution). Therefore,
it is impossible to capture complementarity between alternatives
through a nested logit model. Next, we show that we can capture the
substitutability/complementarity of alternatives through our
welfare-based framework.

In the following discussion, we only consider choice models
$\bq(\bmu)$ that are derived from differentiable choice welfare
functions $w(\bmu)$. We study necessary and sufficient conditions on
the model parameters for a choice model to be substitutable. We
first review the concepts of supermodularity and submodularity:

\begin{definition}[Supermodularity and Submodularity]\label{def:supersubmodular}
A function $f: \mathcal{R}^{n}\mapsto \mathcal{R} \cup \{\infty\}$
is called supermodular if for any $\bx, \by \in \mathcal{R}^{n}$,
$f(\bx \vee \by)+f(\bx \wedge \by)\geq f(\bx)+f(\by)$, where $\bx
\vee \by$ and $\bx \wedge \by$ denote the componentwise maximum and
minimum of $\bx$ and $\by$, respectively. A function $f:
\mathcal{R}^{n}\mapsto \mathcal{R} \cup \{-\infty\}$ is called
submodular if $-f$ is supermodular.
\end{definition}

We have the following theorem:
\begin{theorem}\label{thm:supsubmodular}
Consider a choice model $\bq(\bmu): \mathcal{R}^n \mapsto
\Delta_{n-1}$ that is derived from a differentiable choice welfare
function $w(\bmu)$. Then
\begin{enumerate}
\item $\bq(\bmu)$ is a substitutable choice model if and only if
$w(\bmu)$ is submodular.
\item If $\bq(\bmu)$ is a substitutable
choice model, then there exists an essentially strictly convex
$V(\cdot)$ with $\bar{V}_i(\cdot)$ supermodular on ${\mathcal
R}^{n-1}$ for all $i$, such that
\begin{eqnarray*} \bq(\bmu) =
\arg\max\left\{\bmu^T\bx - V(\bx)
\Big{|}\bx\in\Delta_{n-1}\right\},
\end{eqnarray*}
where
\begin{eqnarray}\label{barv}
\bar{V}_i (\bz) = \left\{\begin{array}{ll}
                                    V\left(z_1,z_2,...,z_{i-1},1-\sum_{j=1}^{n-1} z_j, z_i,...,z_{n-1} \right), & \  \mathrm{if} \  \mb{e}^T \bz \leq 1 \ \mathrm{and} \  \bz \geq 0,  \\
                                    +\infty, & \mathrm{otherwise}.
                                  \end{array}
         \right.
\end{eqnarray}
\end{enumerate}
Furthermore, the reverse is true if $n = 3$.
\end{theorem}

We present the proof of Theorem \ref{thm:supsubmodular} in the
Appendix. Theorem \ref{thm:supsubmodular} provides some sufficient
and necessary conditions for $\bq(\bmu)$ to be substitutable. We
note that the supermodularity of $\bar{V}_i$ has nothing to do with
the supermodularity of $V$. In fact, since $V(\mb{x})$ is only
meaningful on $\Delta_{n-1}$, it can always be modified to be
supermodular by defining $V(\mb{x}) = +\infty$ for all $\mb{x}
\notin \Delta_{n-1}$. The definition of $\bar{V}_i(\cdot)$ reduces a
redundant variable in $V$, making the operations ``$\bx \vee \by$''
and ``$\bx\wedge\by$'' meaningful.

Next we provide an easy-to-check sufficient condition for
a substitutable choice model. We note that in the MDM and the MMM
introduced in Propositions \ref{prop:MDM} and \ref{prop:MMM}, the
corresponding $V(\cdot)$s are separable. The following theorem shows
that the choice models derived from such $V(\cdot)$s are always
substitutable:
\begin{theorem}\label{thm:separable}
If $V(\bx) = \sum_{i \in \mathcal{N}} V_i(x_i)$ on $\Delta_{n-1}$
where $V_i(x_i): [0,1] \mapsto \mathcal{R}^n$ is a strictly convex
function for all $i \in \mathcal{N}$.  Then $\bq(\bmu)$ defined by
\begin{eqnarray}\label{eq:prop:separable}
\bq(\bmu) = \arg\max\left\{\bmu^T\bx - V(\bx) \big{|} \bx \in
\Delta_{n-1}\right\},
\end{eqnarray}
is a substitutable choice model.
\end{theorem}

Another possible choice of  $V(\cdot)$  is a quadratic
function. In that case, we have the following results:
\begin{theorem}\label{thm:V_property}
Consider $\bq(\bmu) = \arg\max\left\{\bmu^T\bx - V(\bx) \big{|} \bx
\in \Delta_{n-1}\right\},$
where $V(\bx) = \bx^T \mb{A} \bx$ is strictly convex with $\mb{A} \succ \mb{0}$. %satisfying $A_{i,j}=a_{i,j}$, for all $i, j \in \mathcal{N}$.
Then $\bar{V}_{i}(\bz)$ for all $i\in \mathcal{N}$ are supermodular
if and only if $A_{jk}-A_{ik}-A_{ij}+A_{ii}\geq 0$ for all distinct
$i,j,k \in \mathcal{N}$, where $A_{ij}$ is the $(i,j)$-th entry of
$A$.
\end{theorem}

Combining Theorems \ref{thm:supsubmodular} and \ref{thm:V_property},
we know that when $n = 3$ and $V(\bx) = \bx^T \mb{A} \bx$ with $\mb{A}\succ \mb{0}$, the
choice model defined by  $\bq(\bmu) = \arg\max\left\{\bmu^T\bx -
V(\bx) \big{|} \bx \in \Delta_{n-1}\right\}$ is substitutable if and
only if
\begin{eqnarray*}
A_{12} + A_{33} \ge A_{13} + A_{23},\quad A_{13} + A_{22} \ge A_{12}
+ A_{23} \mbox{  and  } A_{23} + A_{11} \ge A_{12} + A_{13}.
\end{eqnarray*}
Note that the above condition is different from $\mb{A}$ being
positive semidefinite. Indeed, the following example shows a case
where the choice model is not substitutable even if $V(\bx)$ is
strictly convex and supermodular:
\begin{example}
\label{example:2} Consider $\bq(\bmu) = \arg\max\left\{\bmu^T\bx -
V(\bx) \big{|} \bx \in \Delta_{n-1}\right\},$ where $V(\bx) = \bx^T
\mb{A} \bx$ with
\begin{equation*} \mb{A}=\left(
      \begin{array}{ccc}
        3 & 2 & 0 \\
        2 & 3 & 2 \\
        0 & 2 & 3 \\
      \end{array}
    \right) \succ \mb{0}.
\end{equation*}
It is easy to see that $V(\bx)$ is strictly convex and supermodular.
However, it doesn't satisfy that $A_{13} + A_{22} \ge A_{12} +
A_{23}$. By some further calculations, we obtain that
\begin{eqnarray*}
\bar{V}_2 (\bz) = \bz^T \left(\begin{array}{cc}
               2  & -1 \\
               -1  & 2 \\
              \end{array}\right) \bz - [-2; -2]^T \bz + 3, %, \ \ \  \forall i=1,2,3,$$ % z_1 + z_2 \leq 1, z_1, z_2 \geq 0,$$
\end{eqnarray*}
which is not supermodular.

Therefore $\bq(\bmu)$ is not a substitutable choice model by Theorem
\ref{thm:supsubmodular}. In fact, when we fix $\mu_2 = \mu_3 = 0$
and plot the choice probabilities against $\mu_{1}$ in the range of
values $[-2,2]$ as shown in Figure \ref{fig:ex2}, it is observed
that $q_{3}$ increases in $\mu_{1}$ in the range of $[-1.5,-1]$,
i.e., alternative 3 is complementary to alternative 1 in that
interval. $\hfill\Box$
\begin{figure}[h]
   \centering
    \includegraphics[width = 3.0in]{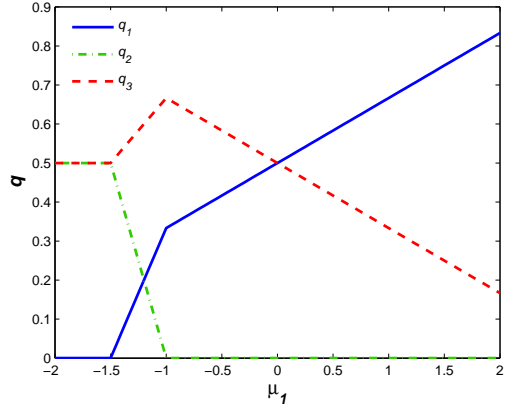}
   \caption{Choice Probabilities in Example \ref{example:2} with $\mu_2 = \mu_3 = 0$}
  \label{fig:ex2}
\end{figure}
\end{example}

In addition to the quadratic example above, we can also easily
generate a non-substitutable choice model through a proper choice
welfare function $w(\cdot)$:
\begin{example}\label{example:1}
Consider the following function:
\begin{equation*}
w(\bmu)=\mathrm{log}\left(e^{\mu_{1}}+e^{\mu_{2}}+e^{\mu_{3}}+e^{0.5(\mu_{1}+\mu_{2})}\right).
\end{equation*}
It is easy to see that $w(\bmu)$ is monotone, translation invariant
and convex, therefore it is a choice welfare function. Also, it is
differentiable. The corresponding choice probability is:
\begin{eqnarray*}
\bq(\bmu)=\frac{1}{e^{\mu_{1}}+e^{\mu_{2}}+e^{\mu_{3}}+e^{0.5(\mu_{1}+\mu_{2})}}\left(e^{\mu_{1}}+\frac{1}{2}e^{0.5(\mu_{1}+\mu_{2})},
e^{\mu_{2}}+\frac{1}{2}e^{0.5(\mu_{1}+\mu_{2})}, e^{\mu_{3}}\right).
\end{eqnarray*}
Furthermore, the second-order derivative of $w(\bmu)$ with respect
to $\mu_{1}$ and $\mu_{2}$ is
\begin{equation*}
\frac{\partial^{2} w(\bmu)}{\partial \mu_{1} \partial
\mu_{2}}=\frac{\partial q_1(\bmu)}{\partial \mu_{2}}=\frac{\partial
q_2(\bmu)}{\partial \mu_{1}} = \frac{
e^{0.5(\mu_{1}+\mu_{2})}(-e^{\mu_{1}}-e^{\mu_{2}}+e^{\mu_{3}}-4e^{0.5(\mu_{1}+\mu_{2})})}{4(e^{\mu_{1}}+e^{\mu_{2}}+e^{\mu_{3}}+e^{0.5(\mu_{1}+\mu_{2})})^{2}}.
\end{equation*}
It is positive if and only if $e^{\mu_{3}} \geq  4  e^{0.5\mu_{1} +
0.5\mu_{2} } +e^{\mu_{1}} + e^{\mu_{2}}$. Therefore, under this
choice model, when $\mu_3$ is large enough (compared to $\mu_1$ and
$\mu_2$), then alternatives $1$ and $2$ will exhibit
complementarity. On the other hand, if $\mu_1$ or $\mu_2$ (or both)
are comparable to $\mu_3$, then they will exhibit substitutability.
Now we give a plausible explanation for this model. Suppose as in
Example 1, we explain $\bmu$ as the number of positive reviews for
each product. Then the above substitutability pattern could be
reasonable if alternatives $1$ and $2$ are both from some relatively
unknown brand (which has very few positive reviews in the history),
while alternative $3$ is from a well-known brand (which in contrast,
has a lot of positive reviews). Then a few more positive reviews on
either alternative $1$ or $2$ will be likely to positively impact
the purchase probability of the other one, since it increases the
overall attractiveness of this brand. On the other hand, if
alternatives $1$ and $2$ have already gained enough positive reviews
in the past, then further increasing the number of positive reviews
of one of them will be more likely to attract the demand from the
other one, rather than from alternative $3$.

To numerically illustrate the above model, we fix $\mu_2=0$, $\mu_3=3$
and plot the choice probability of alternative $2$ as a function of
$\mu_{1}$ in the range of $[-10, 5]$ in Figure \ref{fig:ex1}. From
Figure \ref{fig:ex1}, we can see that alternative $2$ is
complementary to alternative $1$ when $\mu_{1}\in [-10,2]$.
\begin{figure}[h]
   \centering
    \includegraphics[width = 3.0in]{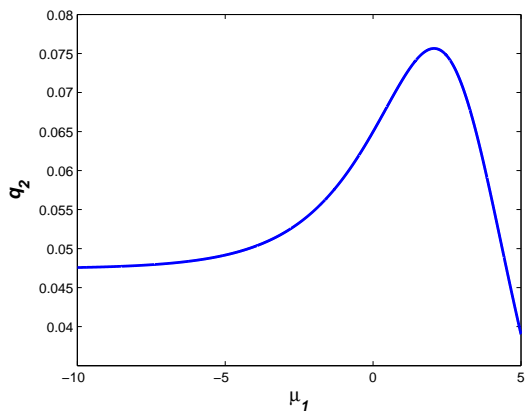}
   \caption{$q_2$ as a Function of $\mu_1$ in Example \ref{example:1}}
  \label{fig:ex1}
\end{figure}
\end{example}
%Next example shows that the supermodularity of $V(\bx)$ itself
%doesn't guarantee the supermodularity of $\bar{V}_{i}(\bz)$ for all
%$i\in \mathcal{N}$. The absence of supermodularity for  some
%$\bar{V}_{i}(\bz)$ when $n=3$ further implies that $\bq_{\mu}$ is
%not a substitutable choice model by Theorem \ref{thm:supsubmodular}.

The above two examples show that by using the welfare-based
framework, it is possible to construct choice models with more
versatile substitution patterns. In addition, the above two examples
further verify that even when $n = 3$, we can construct choice
models that do not have a random utility representation (remember
all random utility models are substitutable choice models).
Therefore, the  welfare-based choice model (thus also the
representative agent model and the semi-parametric models) strictly
subsumes the  random utility model, even for $n = 3$. This
result is an extension of the result obtained by \citet{hofbauer},
which only showed the result for $n \ge 4$.

\section{Constructing New Choice Models from Existing Ones}
\label{sec:construction}

In this section, we show that by using the welfare-based choice
model framework, one can easily construct new choice models from
existing ones. In particular, we provide three transformations below
by which new choice models can be derived from existing ones. In the
following discussions, we use $\bar{\bq}(\cdot)$ to denote existing
welfare-based choice models with choice welfare function
$\bar{w}(\cdot)$, and use $\bq(\cdot)$ and $w(\cdot)$ to denote the
choice probability and the choice welfare function of the new model,
respectively.
%
%The welfare-based choice model exhibits simple structural
%relationship between the utilities of alternatives and the resulting
%choice probability for a sufficiently rich class of choice models.
%In spite of the analytical simplicity of these models, the
%connection between different welfare-based choice models and their
%behavioral properties is not clear so far. In the remainder of this
%section, we provide three transformations by following which new
%welfare-based choice models can be derived from known welfare-based
%choice models. The behavioral properties as well as economic
%intuition underlying the new generated choice model is analyzed.
%More extensive exploration of the behavioral properties of
%welfare-based choice models is provided in Section
%\ref{sec:substitutability}.
%
%The first two transformations (scaled and mixed) are closed under
%RUM, while the third (crossed) transformation is not closed under
%RUM. Thus one can use the crossed transformation to provide new
%choice models given any RUM model. In the following discussion, we
%use $\bar{\bq}(\by)$ to denote a known welfare-based choice model
%with choice welfare function $\bar{w}(\bx)$. Let ${\bq}(\bmu)$
%denote the new welfare-based choice model with choice welfare
%function ${w}(\bmu)$,  which is derived from $\bar{w}(\by)$.
\begin{enumerate}
\item \textbf{Scaling}. Given any existing choice welfare function $\bar{w}(\cdot)$ and any
$\eta > 0$, one can easily verify that $w(\bmu) = \eta
\bar{w}\left({\bmu/\eta}\right)$ is still a choice welfare function.
The corresponding choice model is ${\bq}(\bmu) =
\bar{\bq}\left({\bmu / \eta}\right)$. We note that if
$\bar{\bq}(\cdot)$ has an RUM representation, then $\bq(\cdot)$ also
has an RUM representation with $\mb{\epsilon} = \eta
\mb{\bar{\epsilon}} $. As $\eta$ becomes larger, the difference
among ${\bmu / \eta}$ is smaller and the choice will be more evenly
distributed. As pointed out in \cite{gigerenzer2002bounded}, such
$\eta$ can be used to model the level of rationality of the
individual.

\item \textbf{Mixing}. Let ${B_k, k=1,...,m}$ be a cover of $\mathcal{N}$, i.e., $\cup_k
B_k = \mathcal{N}$. Let $\bar{w}_{k}(\bmu_k)$ be the choice welfare
function on alternatives in $B_k$, with choice probabilities
$\bar{\bq}^{k}(\bmu_k)$. Define
\begin{eqnarray*}
w(\bmu)=\sum_{k=1}^{m}\lambda_{k}\bar{w}_{k}(\bmu_k),
\end{eqnarray*}
where $\lambda_{k}\geq 0$, $\sum_{k=1}^{m}\lambda_{k}=1$. We can
verify that $w(\bmu)$ is a choice welfare function and its
corresponding choice model is
\begin{eqnarray*}
q_i(\bmu)=\sum_{k: i \in B_k} \lambda_{k}
\bar{q}_i^{k}(\bmu_k), \quad \forall i \in \mathcal{N}.
\end{eqnarray*}
If $\bar{\bq}^k(\cdot)$ has an RUM representation for all $k$, then
$\bq(\cdot)$ also has an RUM representation by assuming
$\mb{\epsilon}$ has a mixed distribution of $\mb{\epsilon}_k$,
each with probability $\lambda_k$. This model can be used to model
choice scenarios where there are different segments of customers.
Customers of different segments may only care about a subset of the
products and choose according to a certain choice model. Then the
mixed model is the choice model for the entire population.

\item \textbf{Crossing}. Let $ \mb{A}$ be an $m\times n$ matrix with $A_{ij} \geq 0$ and
$\mb{A} \mb{e}^n = \mb{e}^m$, where $\mb{e}^\ell$ refers to an
$\ell$-dimensional column vector of ones. Given an existing choice
welfare function $\bar{w}(\cdot)$ and its choice probabilities
$\bar{\bq}(\cdot)$, we can easily verify that
\begin{eqnarray*}
w(\bmu)=\bar{w}(\mb{A} \bmu)
\end{eqnarray*}
is still a choice welfare function and the corresponding
welfare-based choice model is
\begin{eqnarray*} \bq(\bmu)=\nabla_{\mbs{\mu}}
w(\bmu)=\mb{A}^{T}\nabla
\bar{w}(\mb{A}\bmu)=\mb{A}^{T}\bar{\bq}(\mb{A}\bmu).
\end{eqnarray*}

An example of such a transformation was in fact shown in Example
\ref{example:1}, where $\bar{w}(\bmu)$ is an MNL model for $4$
alternatives with $\eta = 1$ and
\begin{eqnarray*}
A = \left[\begin{array}{cccc} 1 & 0 & 0 \\
0 & 1 & 0 \\
0 & 0 & 1 \\
0.5 & 0.5 & 0
\end{array}
\right].
\end{eqnarray*}
In light of this example, we note that RUM is not closed under
cross-transformation, i.e., even if $\bar{\bq}(\cdot)$ has an RUM
representation, $\bq(\bmu)$ may not. Thus, the cross-transformation
provides us a way of generating new choice models.
\end{enumerate}

\section{Conclusion}
\label{sec:conclusion}

In this paper, we proposed a new framework for discrete choice
models: the welfare-based choice model, which is based on the idea
of considering the expected utility an individual can get when
facing a set of alternatives. We showed that the welfare-based
choice model is equivalent to the representative agent model and the
semi-parametric  model, thus establishing the equivalence between
the latter two. We also showed that the welfare-based choice model
subsumes the random utility model by relaxing its requirement on
 properties of higher-order cross-partial derivatives
 of the choice welfare function. In particular, we showed that when there
are only two alternatives, the welfare-based choice model is
equivalent to the random utility model. We defined a new concept for
choice models -- substitutability and complementarity -- and showed
that under the new framework, we can construct choice
models with complementary alternatives, thus enabling us to capture new
choice patterns. We believe that this framework is useful for
future studies of choice models.

\bibliographystyle{ormsv080}
\bibliography{choice}

\newpage
\section*{Appendix}
{\noindent\bf Proof of Theorem \ref{thm:choicemodel_equivalence}:} %\\
The equivalence between 1 and 3 directly follows from Theorem 1.
Next we show that $1 \Rightarrow 2$. If $w(\bmu)$ is a
differentiable choice welfare function, by Theorem
\ref{thm:utility_equivalence}, we know that
$$w(\bmu) =
\max \left\{\bmu^T\bx - V(\bx) \big{|} \bx\in \Delta_{n-1}\right\},$$
where
$V(\bx) =  \sup_{\by } \ \{\by^T \bx - w(\by)\}$.
Therefore, $V(\bx)$ is the convex conjugate of $w(\bmu)$. By Theorem
6.3 in \citet{rockafellar}, we know that $w(\bmu)$ is essentially
differentiable if and only if $V(\bx)$ is essentially strictly
convex. Also, from the envelope theorem (see \citealt{mas1995microeconomic}),
$$\nabla w(\bmu) = \nabla_{\mbs{\mu}} \left(\bmu^T\bx - V(\bx)\right)\big{|}_{\mbs{x} = \mbs{x}^*}= \bx^*,$$
where $\bx^* = \argmax \left\{\bmu^T\bx - V(\bx) \big{|} \bx\in
\Delta_{n-1}\right\}$. Therefore,
$$\bq(\bmu) = \nabla w(\bmu) = \argmax \left\{\bmu^T\bx - V(\bx) \big{|} \bx\in \Delta_{n-1}\right\}. $$

Last, we show that $2 \Rightarrow 1 $. Given an essentially strictly
convex $V(\bx)$, by Theorem \ref{thm:utility_equivalence}, we know that
$$w(\bmu) =
\max \left\{\bmu^T\bx - V(\bx) \big{|} \bx\in \Delta_{n-1}\right\}$$
is a choice welfare function. Again, by Theorem 6.3 in
\citet{rockafellar}, we know that $w(\bmu)$ is essentially
differentiable. Moreover, in our case, $w(\bmu)$ is a convex and
finitely valued function in ${\mathcal R}^n$, thus essentially
differentiability is equivalent to differentiability. Again, by
applying the envelope theorem, $\bq(\bmu) = \nabla w(\bmu)$.
Therefore the theorem is proved. $\hfill\Box$\\

%We start with showing that $w(\bmu)$
%is differentiable if and only if $V(\bx)$ is essentially strictly
%convex. By Theorem \ref{thm:utility_equivalence}, we know that
%$V(\bx) =  \sup_{\by } \ \{\by^T \bx - w(\by)\}$ and $w(\bmu) =
%\max \left\{\bmu^T\bx - V(\bx) \big{|} \bx\in \Delta_{n-1}\right\}$.
%Therefore, $V(\bx)$ is the convex conjugate of $w(\bmu)$. By Theorem
%6.3 in \citet{rockafellar}, we know that $w(\bmu)$ is essentially
%differentiable if and only if $V(\bx)$ is essentially strictly
%convex. Moreover, in our case, $w(\bmu)$ is a convex and finitely
%valued function in ${\mathcal R}^n$, thus essentially
%differentiability is equivalent to differentiability. Based on Theorem
%\ref{thm:utility_equivalence}, it remains to show that
%$\hfill \Box$\\
%To prove the equivalence between 1 and 2,

\noindent{\bf Proof of Theorem \ref{thm:boundequivalence}:} First we show the equivalence between 1 and 2. Based on Theorem
\ref{thm:choicemodel_equivalence}, it suffices to prove that $w(\bmu)$
is  superlinear if and only if $V(\bx)$ defined by $\max_{\mbs{y}}\{\by^{T}\bx-w(\by)\}$ is upper bounded. If $w(\bmu)$ is  superlinear,
we have, for any $\bx \in \Delta_{n-1}$,
$$w(\bmu) \geq \sum_{i \in \mathcal{N}} x_i (\mu_i + b_i) = \bx^T \bmu + \bx^T \bb \geq \bx^T \bmu + \min_i \{b_i\} . $$
By reorganizing terms, we have
$$ \bx^T \bmu - w(\bmu) \leq  -\min_i \{b_i\} = \max_i \{-b_i\}. $$
Therefore, $V(\bx) =  \max_{\by } \ \{\by^T \bx - w(\by)\} \leq \max_i \ \{- b_i\}$, i.e., $V(\bx)$ is upper bounded.

To show the other direction, if $V(\bx)$ is upper bounded by a constant $u$, then we have
$$w(\bmu) \geq \max \left\{\bmu^T\bx - u \big{|} \bx\in \Delta_{n-1}\right\} \geq  \mu_i - u , \quad \forall i,$$
i.e., $w(\bmu)$ is superlinear. Therefore, the equivalence between 1 and 2 is proved.

Next we show the equivalence between 1 and 3. We first show that for any  superlinear differentiable choice welfare function $w(\bmu)$, we can find a distribution set $\Theta$ consisting of only distributions with finite expectation such that $w(\bmu)$ can be represented as $w(\bmu)= \sup_{\theta\in \Theta} \mbox{
 } {\mathbb E}_\theta \left[\max_{i\in{\mathcal N}} \mbox{
 }\mu_i + \epsilon_i\right].$

First, since $w(\bmu)$ is convex with $\bq(\bmu) = \nabla w(\bmu)$, we have
\begin{eqnarray}\label{eq:maxofgrad}
w(\bmu) = \sup_{\bz} & \{\bmu^T \bq(\bz) + l(\bz)\},
\end{eqnarray}
where $l(\bz) = w(\bz) - \bz^T\bq(\bz)$. Now we define a
distribution set $\Theta$ that is slightly different from that of
Theorem \ref{thm:utility_equivalence}. Specifically, let $\Theta
=\{\theta_{\mbs{z}} \big{|} \bz\in{\mathcal R}^n\}$, where
$\theta_{\mbs{z}}$ is an $n$-point distribution with ${\mathbb
P}_{\theta_{\mbs{z}}}\left(\bepsilon = \bepsilon_{\mbs{z}}^i\right)
= q_i(\bz)$, $\forall i \in \mathcal{N}$. (Note that by the
monotonicity and the translation invariance properties, $\bq(\bz) =
\nabla w(\bz)$ must satisfy $\bq(\bz) \geq \mathbf{0}$ and $\mb{e}^T \bq(\bz)
=1$.) Here,
\begin{eqnarray*}
\epsilon_{\mbs{z}}^i(j) = \left\{\begin{array}{ll} l(\bz) & \mbox{if }
j = i\\
l(\bz) - M(\bz)& \mbox{if } j\neq i.
\end{array}
\right.
\end{eqnarray*}
where
\begin{equation}\label{eq:Mz}
  M(\bz) = \max \left\{ 1 + \max_{i,j} \ \{z_i - z_j\}, \ {l(\mb{z}) - \min_i \ \{b_i\} \over t^*(\bz)}  \right\},
\end{equation}
with
\begin{equation}\label{eq:tstar}
  t^*(\bz) = \min\{q_i(\mb{z}) | q_i(\mb{z}) >0\}.
\end{equation}
Since $M(\bz) >  z_i - z_j$, for all $i,j$, we have
 $ i = \argmax_{j} \ \{ z_j + \epsilon_{\mbs{z}}^i(j) \}.$
Therefore,
\begin{eqnarray*}
  \mathbb{E}_{\theta_{\mbs{z}}} [ \max_{j} \ \ z_j + \epsilon_j ] = \sum_{i=1}^n q_i(\bz) (z_i + l(\bz))  =
\bz^T \bq(\bz) + l(\bz) = w(\bz).
\end{eqnarray*}
Next we show that:
$$\displaystyle \mathbb{E}_{\theta_{\mbs{z}}} [ \max_{i} \ \ \mu_i + \epsilon_i  ] \leq w(\mb{\mu}), \ \ \forall \mb{\mu}. $$
%Then by (\ref{proofstep1}), we have
%\begin{eqnarray*}
%w(\bmu) = \sup_{\bz} \mbox{ }\bmu^T\bd(\bz) + l(\bz) = \sup_{\bz}
%\mbox{ }{\mathbb E}_{\theta_{\mbs{z}}}[\max_i \ \left\{\mu_i +
%\epsilon_i\right\}] = \sup_{\theta\in \Theta}\mbox{ } {\mathbb
%E}_{\theta} [\max_i \ \left\{\mu_i + \epsilon_i\right\}].
%\end{eqnarray*}
%From
% $$\sup_{\theta \in \Theta} \   \mathbb{E}_\theta [ \max_{i} \ \ \mu_i + \tilde{\epsilon}_i  ] = w(\mb{\mu}).$$
%
%Let's consider the distribution $\theta_{\mbs{\mu},M^*(\mu)}$.
For any given $\bmu$, define
$ k(i) \triangleq  \argmax_{j} \{ \mu_j + \epsilon_{\mbs{z}}^i(j) \}$ (we break ties arbitrarily).
There are two cases:
\begin{enumerate}
  \item For all $i$ such that $q_i(\bz) >0$, $k(i) = i$. In this case, we have
\begin{eqnarray*}
  \mathbb{E}_{\theta_{\mbs{z}}} [ \max_{j} \ \ \{\mu_j + \epsilon_j\} ] = \sum_{i \in \mathcal{N}} q_i(\bz) (\mu_i + l(\bz))  =
\bmu^T \bq(\bz) + l(\bz)
\leq w(\bmu),
\end{eqnarray*}
in which the last inequality is because of the convexity of $w(\cdot)$.
\item There exists some $i$ such that $q_i(\bz) > 0$, but $ k(i) \neq i.$ In this case, from the construction of $\theta_{\mbs{z}}$, we have
  \begin{eqnarray*}
  \mathbb{E}_{\theta_{\mbs{z}}} [ \max_{j} \ \ \mu_j + \epsilon_j ] &=& \sum_{i \in \mathcal{N}, q_{i}(\mbs{z}) >0  } q_{i}(\bz) (\mu_{k(i)} + l(\bz) - M(\bz) {\mathbb I}_{\{k(i) \neq i\}})  \\
  &\leq&  \max_{i} \{\mu_i\} + l(\bz) - t^*(\bz) M(\bz) \\
  &\leq& \max_{i } \{\mu_i\} +  \min_j  \{b_j\} \\
  &\leq& \max_{i} \  \{\mu_i +   b_i\} \\
  &\leq& w(\bmu),
\end{eqnarray*}
where the first inequality follows from the fact that $M(\bz) >0$ and
$\sum_{i\in {\mathcal N}} q_i(\bz) {\mathbb I}_{\{q_i(\mbs{z}) >0, k(i)\neq i\}} \ge t^*(\bz)$,  the second inequality is because of the definition of $M(\bz)$ and the last inequality follows from the definition of superlinear function.
\end{enumerate}
Based on the analysis of these two cases,  we have
$$\displaystyle \mathbb{E}_{ \bepsilon \sim \theta_{\mbs{z}}} [ \max_{i} \ \ \mu_i + \epsilon_i   ] \leq w(\mb{\mu}), \ \ \forall \mb{\mu}. $$
Then by equation (\ref{eq:maxofgrad}) we have
\begin{eqnarray*}
w(\bmu) = \sup_{\bz} \mbox{ } \{\bmu^T\bq(\bz) + l(\bz)\} = \sup_{\bz}
\mbox{ }{\mathbb E}_{\theta_{\mbs{z}}}[\max_i \  \mu_i +
\epsilon_i ] = \sup_{\theta\in \Theta}\mbox{ } {\mathbb
E}_{\theta} [\max_i \  \mu_i + \epsilon_i ].
\end{eqnarray*}
Therefore, we have proved that statement 1 implies statement 3.

Finally, we prove that statement 3 implies statement 1. Suppose there
exists a distribution $\hat{\theta} \in \Theta$ such that
$\mathbb{E}_{\hat{\theta}} |\epsilon_i|  <  +\infty$ for $\forall  i
\in \mathcal{N},$ then for $\bmu \in \mathcal{R}^n$ we have
\begin{eqnarray*}
\sup_{\theta\in \Theta} \mbox{
 } {\mathbb E}_\theta \left[\max_{i\in{\mathcal N}} \mbox{
 }\mu_i + \epsilon_i\right] & \geq & {\mathbb E}_{\hat{\theta}} \left[\max_{i\in{\mathcal N}} \mbox{
 }\mu_i + \epsilon_i\right]
 =  {\mathbb E}_{\hat{\theta}} \left[\mu_j + \epsilon_j\right] = \mu_j + \mathbb{E}_{\hat{\theta}} [\epsilon_j], \quad \forall j.
\end{eqnarray*}
Therefore we can conclude that
$w(\bmu)= \sup_{\theta\in \Theta} \mbox{
 } {\mathbb E}_\theta \left[\max_{i\in{\mathcal N}} \mbox{
 }\mu_i + \epsilon_i\right]$
is  superlinear.

It remains to prove the last statement. We show that for any  $$\bx \in \Delta_{n-1}^\circ \triangleq \{\bx \left| \mb{e}^T \bx =1, x_i > 0, \forall i \in \mathcal{N}\right.\},$$ there exists $\bmu_{\mbs{x}}$ such that $\bq(\bmu_{\mbs{x}})=\nabla w(\bmu_{\mbs{x}})=\bx$. Fix $\bx \in \Delta_{n-1}^\circ$, we consider
\begin{equation}\label{eq:proofthm3no1}
  V(\bx) = \max_{\mbs{\mu}} \ \{\bmu^T \bx - w (\bmu)\}.
\end{equation}
Clearly, $V(\bx) \geq - w (\mb{0})$, since $\bmu=\mb{0}$ is a feasible solution. Moreover, since $w (\bmu)$ is translation invariant, we can restrict the feasible region of (\ref{eq:proofthm3no1}) to  $\mathcal{L} \triangleq\{\bmu | \mb{e}^T \bmu=0\}$. For all $\bmu \in \mathcal{L}$, we have $\mu_j \leq 0$ for some $j \in \mathcal{N}$. Thus
$$\bmu^T\bx \leq \sum_{i \neq j} \mu_i x_i  \leq \sum_{i \neq j} x_i  \max_{k} \{\mu_k\} \leq (1-\min_i \{x_i\}) \max_{k} \{\mu_k\}.$$
However, by superlinearity of $w(\bmu)$, we have:
$$w(\bmu) \geq \max_{k} \{\mu_k + b_k\} \geq \max_{k} \{\mu_k\} + \min_{k}\{b_k\}.$$
Thus, for all $\bmu \in \mathcal{L}$, we have:
$$\bmu^T \bx - w (\bmu) \leq -\min_i \{x_i\} \max_{k} \{\mu_k\} - \min_{k}\{b_k\} .$$
Let $K= {w(\mbs{0}) - \min_{k}\{b_k\} \over \min_i \{x_i\}}$. In
order for $\bmu$ to be optimal to (\ref{eq:proofthm3no1}), by the
above arguments, we would have $\mu_i \leq K$ for all $i$. Thus we
can further restrict the feasible set of (\ref{eq:proofthm3no1}) to
$\{\bmu | \mb{e}^T \bmu=0, \mu_i \leq K \ \forall i \in \mathcal{N}
\}$, which is a compact set. Since $w(\bmu)$ is continuous, there
exists $\bmu_{\mbs{x}} \in \{\bmu | \mb{e}^T \bmu=0, \mu_i \leq K \
\forall i \in \mathcal{N} \}$ that attains maximum in problem
(\ref{eq:proofthm3no1}). By the first-order necessary condition,
$\nabla w(\bmu_{\mbs{x}})=\bx. $ This concludes the proof.
 $\hfill\Box$\\

{\noindent\bf Proof of Proposition \ref{prop:basicpropertysub}:}
Since $w(\bmu)$ is convex and differentiable, for any $\bmu \in {\mathcal R}^n$ and any $t > 0$, we have
\begin{equation*}
  \begin{array}{l}
   w(\bmu + t \mb{e_i} ) - w(\bmu)  \geq   t \mb{e_i^T} \nabla w(\bmu) = t q_i(\bmu), \vspace{4pt}\\
    w(\bmu) - w(\bmu + t \mb{e_i} )   \geq   - t \mb{e_i^T} \nabla w (\bmu + t \mb{e_i} ) = - t q_i(\bmu + t \mb{e_i}). \vspace{4pt}
  \end{array}
\end{equation*}
From these two inequalities, we have $ q_i(\bmu + t \mb{e_i}) - q_i(\bmu) \geq 0$, for all $t > 0$ and $\bmu$. Thus, alternative $i$ is complementary to itself.

Furthermore, if $w(\bmu)$ is second-order continuously
differentiable, then we have ${\partial q_i \over \partial \mu_j } =
{\partial^2 w \over \partial \mu_i \partial \mu_j } =  {\partial^2 w
\over \partial \mu_j \partial \mu_i }   = {\partial q_j \over
\partial \mu_i }$. Thus, if alternative $i$ is substitutable
(complementary, resp.) to alternative $j$ at
$\bmu$, then alternative $j$ is substitutable (complementary, resp.) to alternative $i$ at $\bmu$.   $\hfill\Box$\\

%{\noindent\bf Proof of Theorem \ref{thm:randomutilitysub}:} In a
%random utility model, the probability of choosing alternative $i$ is
%$q_i(\bmu) = {\mathbb P}_{\theta}\left(i =\arg\max_{k\in{\mathcal
%N}} \ (\mu_k + \epsilon_k)\right),$ which is decreasing in $\mu_{j}$
%for all $j\neq i$ and $\bmu$. Therefore all pairs of alternatives
%$i$ and $j$ ($i\neq j$) must be substitutable to each other in a
%random utility model. $\hfill\Box$

{\noindent\bf Proof of Theorem \ref{thm:supsubmodular}:}
In this proof, we use the following lemma  from \cite{murota2003discrete}.
\begin{lemma}
\label{lemma:supersub} Let $f: \mathcal{R}^{n}\mapsto \mathcal{R}
\cup \{\infty\}$ be a function such that there exists at least one
$\bmu$ such that $f(\bmu)<\infty$. Let $g(\mb{x}) = \max_{ \bmu }
\quad \{\bmu^T \bx - f (\mb{\mu})\}$ be the \textit{convex
conjugate} of $f$. We have
\begin{enumerate}
  \item If $f$ is submodular, then $g$ is supermodular.
  \item If $n=2$ and $f$ is supermodular, then $g$ is submodular.
\end{enumerate}
%Moreover, if $g$ is differentiable at $\mb{x}$, then
%${\partial^2 g(\mbs{x}) \over \partial x_i \partial x_j } \geq 0, \forall i \neq j$.
\end{lemma}

Now we use this lemma to prove the theorem. To prove the first part, by \cite{simchi2014convexity}, a differentiable function $w(\bmu)$ is submodular in $\bmu$ if and only if ${\partial w(\mbss{\mu}) \over \partial \mu_i }$ is decreasing in $\mu_j$ for all $i \neq j$. By the definition of $\bq(\bmu) = \nabla w(\bmu)$, the result holds.

For the second part, let $V(\mb{x}) = \max_{ \bmu } \ \{\bmu^T \bx - w (\mb{\mu})\}$ be the convex conjugate of $w(\bmu)$. From Theorem \ref{thm:choicemodel_equivalence}, $V(\mb{x})$ is essentially strictly convex and $$\bq(\bmu) = \arg\max\left\{\bmu^T\bx - V(\bx) \Big{|}\bx\in\Delta_{n-1}\right\}.$$
For any $\by \in \mathcal{R}^{n-1}$ and $i\in \mathcal{N}$, define $f_{i}(\by)=w(y_{1},y_{2},...,y_{i-1},0,y_{i},...,y_{n-1})$.  Also define $\mb{{\mu}_{-i}}=(\mu_{1},...,\mu_{i-1},\mu_{i+1},...,\mu_{n})$, then we have
\begin{eqnarray*}
\bar{V}_i (\bz) & = &\max_{\mbs{\mu}} \quad  \{\mb{\mu}_{-i}^T \mb{z} + \mu_{i} (1-\mb{e}^T \mb{z}) - w (\mb{\mu})\}\\
 & = & \max_{\mbs{\mu}, \mu_{i}=0} \quad \{ \mb{\mu}_{-i}^T \mb{z} + \mu_{i} (1-\mb{e}^T \mb{z}) - w (\mb{\mu}) \}\\
 &=&\max_{\mbs{y}} \quad \{\mb{y}^T \mb{z} - f_{i}(\mb{y})\},
\end{eqnarray*}
%$$V(1-\mb{e}^T \mb{x_{-i}},  \mb{x_{-i}}) = \max_{\mbs{\mu_{-i}}}$$
where the second equality is due to the translation invariance property of $w (\mb{\mu})$. The submodularity of $w(\bmu)$ implies the submodularity of $f_{i}(\by)$ for all $i \in \mathcal{N}$. Thus $\bar{V}_i (\bz)$, as the convex conjugate of $f_{i}(\by)$, is supermodular by Lemma \ref{lemma:supersub}.

For the last statement, since $V(\cdot)$ is an essentially strictly
convex function, $\bq(\bmu) = \arg\max\left\{\bmu^T\bx - V(\bx)
\Big{|}\bx\in\Delta_{n-1}\right\}$ is well-defined. By Theorem
\ref{thm:choicemodel_equivalence}, $\bq(\bmu) = \nabla w(\bmu)$
where $w(\bmu) = \sup\left\{\bmu^T\bx - V(\bx) \big{|} \bx\in
\Delta_{n-1}\right\}$. For any $\by \in \mathcal{R}^{n-1}$ and $i\in
\mathcal{N}$, define
$f_{i}(\by)=w(y_{1},y_{2},...,y_{i-1},0,y_{i},...,y_{n-1})$.  Also
define $\mb{{x}_{-i}}=(x_{1},...,x_{i-1},x_{i+1},...,x_{n})$, then
we have
\begin{eqnarray*}
f_{i}(\by)  & = &\max_{\mbs{x} \in \Delta_{n-1}} \quad  \left\{\mb{x}_{-i}^T \mb{y} + 0 (1-\mb{e}^T \mb{x}_{-i}) -  V(\bx) \right\}\\
 & = & \max_{\mbs{x} \in \Delta_{n-1}} \quad \left\{ \mb{x}_{-i}^T \mb{y} -  \bar{V}_i (\mb{x}_{-i}) \right\} \\
  & = & \max_{\mbs{x}} \quad \left\{ \mb{y}^T \mb{x}_{-i}  -  \bar{V}_i (\mb{x}_{-i}) \right\} \\
 &=&\max_{\mbs{z}} \quad \left\{ \mb{y}^T \mb{z} - \bar{V}_i(\mb{z}) \right\},
\end{eqnarray*}
where the third equality holds since $ \bar{V}_i (\mb{x}_{-i}) = +\infty$ for all $\mb{x} \notin \Delta_{n-1}$. From Lemma \ref{lemma:supersub}, given that $n=3$ and thus $\by \in \mathcal{R}^2$, $f_{i}(\by)$ is submodular. It remains to show that $w(\bmu)$ is also submodular. According to Theorem \ref{thm:supsubmodular}, it suffices to show that $q_i(\bmu)$ is locally decreasing with $\mu_j$ for all $j\neq i$ for all $\bmu$. Fix $i,j$ and let $k \neq i, j $. We assume $i>j$ without loss of generality. We have $q_i(\bmu -\mu_k \mb{e} ) = q_i(\bmu) $ from translation invariance property. But $q_i(\bmu -\mu_k \mb{e} ) = {\partial f_{k}(\mu_i-\mu_k, \mu_j - \mu_k ) \over \partial \mu_i} $ is non-decreasing with $\mu_j$ due to the submodularity of $f_{k}$. Thus $w(\bmu)$ is submodular and $\bq(\bmu) = \nabla w(\bmu)$ is a substitutable choice model.
$\hfill\Box$ \\

{\noindent\bf Proof of Theorem \ref{thm:separable}:} We first consider
the case where $V_i(x_i)$ is differentiable for all $i \in
\mathcal{N}$. Let $\lambda(\bmu)$ be the Lagrangian multiplier
of the constraint $\sum_i x_i = 1$. The KKT conditions (see
\citealt{bertsekas}) for problem (\ref{eq:prop:separable}) can be
written as:
$$\begin{array}{cl}
    \mu_i - V_i'(q_i(\bmu)) - \lambda(\bmu) \leq 0, &\ \forall i \in \mathcal{N}; \vspace{4pt} \\
    \mu_i - V_i'(q_i(\bmu)) - \lambda(\bmu) =  0, &\ \forall i \ \mathrm{s.t.} \  q_i(\bmu) \neq 0; \vspace{4pt}  \\
     q_i(\bmu) \geq 0, &\ \forall i \in \mathcal{N}; \vspace{4pt} \\
     \sum_{i \in \mathcal{N}} q_i(\bmu) = 1. \vspace{4pt}
  \end{array}
$$
Now we consider any two points $\bmu_0$ and $\bmu_0 + t \mb{e_i}$ where $\mb{e_i}$ is a unit vector along the $i$-th coordinate axis  and $t > 0$.  Suppose that there exists a $j \neq i$ such that $q_j(\bmu_0 + t \mb{e_i}) > q_j(\bmu_0)$.  Since $V_j$ is strictly convex, $V_j'(q_j(\bmu_0 + t \mb{e_i})) > V_j'(q_j(\bmu_0))$. There are two possible cases for $q_j(\bmu_0)$:
\begin{itemize}
  \item $q_j(\bmu_0)>0$: In  this case, we have $\mu_j - V_j'(q_j(\bmu_0+ t \mb{e_i})) - \lambda(\bmu_0+ t \mb{e_i}) = 0$ and $\mu_j - V_j'(q_j(\bmu_0)) - \lambda(\bmu_0) = 0$, therefore, we have $\lambda(\bmu_0 + t \mb{e_i}) < \lambda(\bmu_0 )$.
  \item $q_j(\bmu_0)=0$:
In this case, $\mu_j - V_j'(q_j(\bmu_0)) - \lambda(\bmu_0) \leq 0$, which implies that $\mu_j - V_j'(q_j(\bmu_0 + t \mb{e_i})) - \lambda(\bmu_0) < 0$. But $\mu_j - V_j'(q_j(\bmu_0 + t \mb{e_i})) - \lambda(\bmu_0 + t \mb{e_i}) = 0$, we have
$\lambda(\bmu_0 + t \mb{e_i}) < \lambda(\bmu_0 )$.
\end{itemize}
In both cases, $\lambda(\bmu_0 + t \mb{e_i}) < \lambda(\bmu_0 )$. This implies that $q_j(\bmu_0 + t \mb{e_i}) \geq q_j(\bmu_0)$ for all $j \neq i$. Note that we also have $q_i(\bmu_0 + t\mb{e_i}) > q_i(\bmu_0)$ by Proposition \ref{prop:basicpropertysub}. Therefore, we have $\sum_{j \in \mathcal{N}} q_j(\bmu_0+ t \mb{e_i}) > \sum_{j \in \mathcal{N}} q_j(\bmu_0) = 1$, which contradicts with that $\bq(\bmu_0 + t\mb{e_i} )\in \Delta_{n-1}$.
 Thus we have $q_j(\bmu_0 + t \mb{e_i}) \leq q_j(\bmu_0)$ for all $j \neq i$. Since this is true for all $\bmu_0$ and $t >0$, $\bq$ is substitutable.

If $V_i(x_i)$ is not differentiable, we need to replace the derivative with the subgradient in the above argument. Since $V_i$ is strictly convex, $g_1>g_2$ for all $g_1 \in \partial V_i(x_1)$ and $g_2 \in \partial V_i(x_2)$ if $x_1 > x_2$,  the above argument is still valid.
$\hfill\Box$\\

{\noindent\bf Proof of Theorem \ref{thm:V_property}:} For $i\in
\mathcal{N}$, $\bar{V}_{i}$ is an $n-1$ variate quadratic function.
Let $H^i$ denote the Hessian matrix of $\bar{V}_{i}(\bz)$. For
$j,k\in \{1,2,...,n-1\}$ and $j\neq k$,  the off-diagonal element
$H^i_{j,k}=A_{\tilde{j},\tilde{k}}-A_{i,\tilde{k}}-A_{i,\tilde{j}}+A_{i,i}$,
where $$\tilde{j}=\left\{\begin{array}{cl}
                   j, \ & \mathrm{if} \ j<i, \\
                   j+1, \ & \mathrm{if} \ j\geq i;
                 \end{array}
\right. \quad \ \mathrm{and} \ \quad \tilde{k}=\left\{\begin{array}{cl}
                   k, \ & \mathrm{if} \ k<i, \\
                   k+1, \ & \mathrm{if} \ k\geq i.
                 \end{array}
\right.$$
Thus, $\bar{V}_{i}(\bz)$ is supermodular if and only if $H^i_{j,k}\geq 0$ for all $j,k\in \{1,2,...,n-1\}$ and $j\neq k$, which is equivalent to $A_{j,k}-A_{i,k}-A_{i,j}+A_{i,i}\geq 0$ for all distinct $i,j,k \in \mathcal{N}$. $\hfill\Box$

\end{document}